\documentclass[11pt]{amsart}

\usepackage{geometry}
\usepackage{todonotes}

\usepackage{amsfonts, amssymb, amscd}
\numberwithin{equation}{section}

\usepackage[symbol]{footmisc}

\usepackage{bm}
\usepackage{verbatim}
\usepackage{mathrsfs}
\usepackage{graphicx}
\usepackage{tikz-cd}
\usepackage{subcaption}
\usepackage{listings}
\usepackage{subfiles}
\usepackage[toc,page]{appendix}
\usepackage{mathtools}
\usepackage{comment}
\usepackage{enumerate}
\usepackage{enumitem}
\usepackage[all]{xy}

\usepackage{graphicx}
\graphicspath{{images/}}

\usepackage{appendix}
\usepackage{hyperref}
\hypersetup{
    colorlinks=true,
    citecolor=red,
    linkcolor=blue,
    filecolor=magenta,      
    urlcolor=red,
}
\newcommand{\bb}{\bm{b}}
\newcommand{\Mm}{{\bf{M}}}

\newcommand{\Cc}{\mathbb{C}}

\newcommand{\Pp}{\mathbb{P}}
\newcommand{\Qq}{\mathbb{Q}}

\newcommand{\Rr}{\mathbb{R}}

\newcommand{\Center}{\operatorname{center}}

\newcommand{\Exc}{\operatorname{Exc}}

\newcommand{\Nklt}{\operatorname{Nklt}}

\newcommand{\Supp}{\operatorname{Supp}}

\newcommand{\mult}{\operatorname{mult}}

\newcommand{\lf}{\lfloor}
\newcommand{\rf}{\rfloor}

\newcommand{\Oo}{\mathcal{O}}
\newcommand{\Ii}{\mathcal{I}}

\newcommand{\Pic}{\mathrm{Pic}}

\newtheorem{thm}{Theorem}[section]

\newtheorem{lem}[thm]{Lemma}
\newtheorem{prop}[thm]{Proposition}

\newtheorem{claim}[thm]{Claim}

\theoremstyle{definition}
\newtheorem{defn}[thm]{Definition}

\theoremstyle{definition}

\newtheorem{ex}[thm]{Example}

\theoremstyle{definition}

\begin{document}

\title{Relative Nakayama-Zariski decomposition and minimal models of generalized pairs}
\author{Jihao Liu and Lingyao Xie}

\address{Department of Mathematics, Northwestern University, 2033 Sheridan Rd, Evanston, IL 60208, USA}
\email{jliu@northwestern.edu}

\address{Department of Mathematics, The University of Utah, Salt Lake City, UT 84112, USA}
\email{lingyao@math.utah.edu}

\subjclass[2020]{14E30,14C20.14E05,14J17}
\date{\today}

\begin{abstract}
We prove some basic properties of the relative Nakayama-Zariski decomposition. We apply them to the study of lc generalized pairs. We prove the existence of log minimal models or Mori fiber spaces for (relative) lc generalized pairs polarized by an ample divisor. This extends a result of Hashizume-Hu to generalized pairs. We also show that, for any lc generalized pair $(X,B+A,\Mm)/Z$ such that $K_X+B+A+\Mm_X\sim_{\Rr,Z}0$ and $B\geq 0,A\geq 0$, $(X,B,\Mm)/Z$ has either a log minimal model or a Mori fiber space. This is an analogue of a result of Birkar/Hacon-Xu and Hashizume in the category of generalized pairs, and is later shown to be crucial to the proof of the existence of lc generalized flips in full generality.
\end{abstract}

\maketitle

\tableofcontents
\section{Introduction}

We work over the field of complex numbers $\mathbb C$.

The theory of \emph{generalized pairs} (\emph{g-pairs} for short) was introduced by Birkar and Zhang in \cite{BZ16} to tackle the effective Iitaka fibration conjecture. The structure of g-pairs naturally appears in the canonical bundle formula and sub-adjunction formulas \cite{Kaw98,FM00}. This theory has been used in an essential way in the proof of the Borisov-Alexeev-Borisov conjecture \cite{Bir19,Bir21a}. We refer the reader to \cite{Bir21b} for a more detailed introduction to the theory of g-pairs.

Recently, there is significant progress towards the minimal model program theory for generalized pairs. In particular, in \cite{HL21a}, Hacon and the first author proved the cone theorem, contraction theorem, and the existence of flips for $\Qq$-factorial lc g-pairs. However, some related results on the termination of flips and the existence of log minimal models and good minimal models for generalized pairs remain unknown. For example, we have the following results in the setting of usual pairs:

\begin{thm}[{\cite[Theorem 1.5]{HH20}}]\label{thm: hh20 1.5}
Let $(X,B)/Z$ be a pair and $A\geq 0$ an ample$/Z$ $\Rr$-divisor such that $(X,\Delta:=B+A)$ is lc and $Z$ is normal quasi-projective. Then $(X,\Delta)/Z$ has a good minimal model or a Mori fiber space.
\end{thm}

\begin{thm}[{\cite[Theorem 1.1]{Has19}; see \cite{Bir12,HX13} for the $\Qq$-coefficient case}]\label{thm: has19 1.1}
Let $(X,B)/Z$ be a pair and $A\geq 0$ an $\Rr$-divisor such that $(X,B+A)$ is lc, $Z$ is normal quasi-projective, and $K_X+B+A\sim_{\Rr,Z}0$. Then:
    \begin{enumerate}
        \item $(X,B)/Z$ has either a Mori fiber space or a log minimal model $(Y,B_Y)/Z$.
        \item If $K_Y+B_Y$ is nef$/Z$, then $K_Y+B_Y$ is semi-ample$/Z$.
        \item If $(X,B)$ is $\Qq$-factorial dlt, then any $(K_X+B)$-MMP$/Z$ with scaling of an ample$/Z$ $\Rr$-divisor terminates.
    \end{enumerate}
\end{thm}

In this paper, we further investigate the minimal model program for generalized pairs. We prove the following results, which can be considered as analogues of Theorems \ref{thm: hh20 1.5} and \ref{thm: has19 1.1} respectively:

\begin{thm}\label{thm: Has22a 3.17 rel ver intro}
Let $(X,B,\Mm)/U$ be an NQC lc g-pair and $A\geq 0$ an ample$/U$ $\Rr$-divisor such that $(X,\Delta:=B+A,\Mm)$ is lc. Then
\begin{enumerate}
    \item $(X,\Delta,\Mm)/U$ has a log minimal model or a Mori fiber space, and
    \item if $\Mm_X$ is $\Rr$-Cartier, then  $(X,\Delta,\Mm)/U$ has a good minimal model or a Mori fiber space.
\end{enumerate}
\end{thm}

\begin{thm}\label{thm: bir12 1.1 gpair}
Let $(X,B,\Mm)/U$ be an NQC lc g-pair such that $X\rightarrow U$ is a projective morphism between normal quasi-projective varieties, and $A\geq 0$ an $\Rr$-divisor such that $(X,B+A,\Mm)$ is lc and $K_X+B+A+\Mm_X\sim_{\Rr,U}0$. Then
\begin{enumerate}
    \item $(X,B,\Mm)/U$ has a log minimal model or a Mori fiber space, and
    \item if $(X,B,\Mm)$ is $\Qq$-factorial dlt, then any $(K_X+B+\Mm_X)$-MMP$/U$ with scaling of an ample$/U$ $\Rr$-divisor terminates.
\end{enumerate}
\end{thm}

Theorems \ref{thm: Has22a 3.17 rel ver intro} and \ref{thm: bir12 1.1 gpair} have played important roles in the minimal model program theory for lc generalized pairs, especially the existence of generalized lc flips. See the \hyperlink{Postscript}{Postscript} for details.

Note that when $\Mm=\bf{0}$, Theorem \ref{thm: Has22a 3.17 rel ver intro} is exactly Theorem \ref{thm: hh20 1.5} and Theorem \ref{thm: bir12 1.1 gpair} is exactly Theorem \ref{thm: has19 1.1}(1)(3). For technical reasons, at the moment, we cannot remove the ``$\Mm_X$ is $\Rr$-Cartier" assumption in Theorem \ref{thm: Has22a 3.17 rel ver intro}(2).

We still expect the analogue of Theorem \ref{thm: has19 1.1}(2) to be true. That is, we expect that any log minimal model of $(X,B,\Mm)/Z$ is a good minimal model (of) a generalized pair $(X,B,\Mm)/Z$ as in Theorem \ref{thm: bir12 1.1 gpair} is in fact good; see the first paragraph of the Postscript. This is because such $K_X+B+\Mm_X$ is log abundant$/U$ with respect to $(X,B,\Mm)$ by Theorem \ref{thm: Has22b 4.1 rel ver} below. However, the following example shows that the question is very subtle as ``log abundance" does not imply semi-ampleness in general for lc g-pairs:

\begin{ex}\label{ex: log abundant not semiample}
Let $C_0$ be a nodal cubic in $\Pp^2$ and $l$ the hyperplane class on $\Pp^2$. Let $P_1,P_2,...,P_{12}$ be twelve distinct points on $C_0$ which are different from the nodal point. Let
$$
\mu:X=\text{Bl}_{\{P_1,...,P_{12}\}}\to\Pp^2
$$
be the blow-up of $\Pp^2$ at the chosen points with the exceptional divisor $E=\sum_{i=1}^{12}E_i$, where $E_i$ is the prime exceptional divisor over $P_i$ for each $i$. Let $H:=\mu^*l$ and $C:=\mu^{-1}_*C_0$. Then $C\cong C_0$, $C\in|3H-E|$, and $K_X+C=\mu^*(K_{\Pp^2}+C_0)=0$. Moreover, we have $C_0^2=9$ and $C^2=-3$.

We consider the big divisor $M=4H-E\sim H+C$. Since $H$ is semi-ample and $M\cdot C=0$, $M$ is nef. Notice that $\Oo_C(M)=\Oo_{C_0}(4l-\sum_{i=1}^{12}P_i)$ and $\Pic^0(C)\cong\mathbb G_m$, where $\mathbb G_m$ is the multiplication group of $\Cc^*$.

\begin{enumerate}
    \item Suppose that $P_1,...,P_{12}$ are in general position so that $\Oo_C(M)$ is a non-torsion line bundle in $\Pic^0(C)$. Then $M$ can never be semi-ample since $M|_C$ is not. However, the normalization $C^n$ of $C$ is $\Pp^1$, so $M|_{C^n}$ is semi-ample. This gives an lc g-pair $(X,C,\Mm:=\overline{M})$ such that both $M$ and $K_X+C+M\sim M$ are nef and log abundant with respect to $(X,C,\Mm)$, but $K_X+C+M$ is not semi-ample. One can further take the blow-up of the nodal point and take the crepant pullback to make each lc center normal.
    \item Suppose that $P_1,...,P_{12}$ are the intersection points of $C_0$ with a general quartic curve $Q_0\in|4l|$. Let $Q$ be the birational transform of $Q_0$ on $X$. Then $M\sim Q\sim H+C$ is semi-ample and defines a projective birational contraction $f: X\to Y$ which contracts exactly the nodal curve $C$. Let $M'=H-3E_1$. Then $M'\cdot C=0$ and $\Oo_C(M')=\Oo_{C_0}(l-3P_1)$ is a non-torsion line bundle since $Q_0$ is general. Therefore $M'$ is not $\Qq$-linearly equivalent to 0 over $Y$ (which also implies that $f(M')$ is not $\Qq$-Cartier). This gives an lc g-pair $(X,C,\Mm':=\overline{M'})/Y$ such that  both $M'$ and $K_X+C+M'\sim M'$ are log abundant and numerically trivial over $Y$ but $K_X+C+M'$ is not semi-ample over $Y$.
\end{enumerate}
\end{ex}

We refer the reader to \cite{BH22} for some other interesting examples on the failure of positivity results for generalized pairs.

\medskip

To prove our main theorems, the central idea is to combine the methods in \cite{Has22a} (some originated in \cite{Has20,Has22b,HH20}) and \cite{HL21a}. In particular, we need to generalize many results in \cite{Has22a} for projective varieties $X$ to normal quasi-projective varieties $X$ equipped with projective morphisms $\pi: X\rightarrow U$. Despite their similarities, a major difficulty is the use of the Nakayama-Zariski decomposition \cite[III. \S 1]{Nak04}, which is usually applied to projective varieties only. It is important to remark that the relative Nakayama-Zariski decomposition \cite[III. \S 4]{Nak04} does not always behave as good as the global Nakayama-Zariski decomposition (see \cite{Les16}), and we lack references for even the most basic properties of them. In this note, we will study the behavior and basic properties of the relative Nakayama-Zariski decomposition. We refer the reader to \cite{LT22b} for further applications of the relative Nakayama-Zariski decomposition on the minimal model theory for generalized pairs.

\medskip

\noindent\textit{Idea of the proof}. It is important to notice that Theorems \ref{thm: Has22a 3.17 rel ver intro} and \ref{thm: bir12 1.1 gpair} both have some ``$\bb$-log abundant" conditions:
\begin{enumerate}
    \item In Theorem \ref{thm: Has22a 3.17 rel ver intro}, possibly replacing $(X,B,\Mm)$ with $(X,B,\Mm+\frac{1}{2}\bar A)$ and $A$ with $\frac{1}{2}A$, we may assume that $\Mm$ is $\bb$-log abundant with respect to $(X,B,\Mm)$.
    \item In Theorem \ref{thm: bir12 1.1 gpair}, $K_X+B+A+\Mm_X$ is automatically $\bb$-log abundant$/Z$ as it is $\Rr$-linearly trivial over $Z$.
\end{enumerate}
Therefore, one important goal of this paper is to study the minimal model program for g-pairs $(X,B,\Mm)$ with $\bb$-log abundant nef part $\Mm$ or with log abundant $K_X+B+\Mm_X$. Despite the technicality, the condition ``$\bb$-log abundant" is actually a very natural condition as it is preserved under adjunction. The key idea to study the minimal model program for such g-pairs is the following:
\begin{itemize}
    \item By applying the Iitaka fibration and the generalized canonical bundle formula, we reduce the questions to the cases when either $\kappa_{\iota}(X/U,K_X+B+\Mm_X)=0$ or $\kappa_{\iota}(X/U,K_X+B+\Mm_X)=\dim X-\dim U$ (see Section 4).
    \item When the invariant Iitaka dimension is $0$, by abundance, the minimal model program behaves well (cf. Lemma \ref{lem: has19 3.2 step 3 abu ver}). So we can reduce the question to the case when $K_X+B+\Mm_X$ is big$/U$.
    \item If $(X,B,\Mm)$ is klt then we can apply \cite[Lemma 4.4(2)]{BZ16}. Otherwise, by induction on dimension, we can apply special termination results near $\Nklt(X,B,\Mm)$.
\end{itemize}

\noindent\textit{Structure of the paper}. In Section 2, we introduce some preliminary results. In particular, we will recall some results on the minimal model program for generalized pairs that are already included in \cite[Version 2, Version 3]{HL21a} (but may not appear in the published version). In Section 3, we study the basic behavior of the relative Nakayama-Zariski decomposition. In Section 4, we use the Iitaka fibration and the generalized canonical bundle formula to simplify the question. In Section 5,6 and 7, we use the relative Nakayama-Zariski decomposition to prove analogues of most results in \cite[Section 3]{Has22a} (Section 5), \cite[Theorem 3.14]{Has22a} (Section 6), and \cite[Theorem 4.1]{Has22b} (Section 7) respectively. In Section 8, we prove Theorems \ref{thm: Has22a 3.17 rel ver intro} and \ref{thm: bir12 1.1 gpair}.

\medskip

\noindent\textbf{Acknowledgement}. The authors would like to thank Christopher D. Hacon, Jingjun Han, Junpeng Jiao, Vladimir Lazi\'c, Yuchen Liu, Yujie Luo, Fanjun Meng, Nikolas Tsakanikas, and Qingyuan Xue for useful discussions. The authors would like to thank Kenta Hashizume for useful comments on an earlier manuscript of this paper. Part of the work is done during the visit of the first author to the University of Utah in March and April 2022, and the first author would like to thank their hospitality. The second author is partially supported by NSF research grants no: DMS-1801851, DMS-1952522 and by a grant from the Simons Foundation; Award Number: 256202.

Some parts of this note has overlap with results in \cite[Version 2 or Version 3]{HL21a}. Since these results are not expected to be published in the final version of \cite{HL21a} due to the simplification of the proofs of the main theorems of \cite{HL21a}, for the reader's convenience, we include some of the results of \cite[Version 2 or Version 3]{HL21a} in this paper and provide detailed proofs. The authors would like to thank Christopher D. Hacon for granting the text overlap.

We thank Xiaowei Jiang for useful comments for the first version of the paper.

We wish to thank the anonymous referee for useful
comments and suggestions that helped the authors improve the clarity of this work.

\medskip

\noindent\textbf{\hypertarget{Postscript}{Postscript}}. After the first version of the paper appeared on the arXiv, the authors proved a stronger version of Theorem \ref{thm: bir12 1.1 gpair} in \cite[Theorem 1.1]{LX22}, which shows that the log minimal model $(Y,B_Y,\Mm)/Z$ is essentially a good minimal model of $(X,B,\Mm)/Z$. This is crucial for the complete solution of the existence of flips for lc generalized pairs \cite[Theorem 1.2]{LX22}, which removes the $\mathbb R$-Cartier assumption of $\Mm_X$ as in \cite[Theorem 1.2]{HL21a}. Although the proof of \cite{LX22} heavily relies on this paper, we decided to write and submit them as two separate papers, as this paper contains most technical results that we need while \cite{LX22} mainly focuses on establishing a Koll\'ar-type gluing theory. 

We also remark that the second author and N. Tsakanikas proved a stronger version of Theorem \ref{thm: Has22a 3.17 rel ver intro}, removing the $\Rr$-Cartierness assumption of $\Mm_X$ in Theorem \ref{thm: Has22a 3.17 rel ver intro}(2), see \cite[Theorem F]{TX23}. The proof of \cite[Theorem F]{TX23} relies on \cite{LX22,Xie22}  which in turn rely on this paper. Therefore, we will avoid citing results from \cite{LX22,Xie22,TX23} in this paper.

\section{Preliminaries}

We adopt the same notation as in \cite{KM98,BCHM10}. For g-pairs, we adopt the same notation as in \cite{HL21a}, which is the same as \cite{FS20,Has22a} except that we use ``$a(E,X,B,\Mm)$ instead of ``$a(E,X,B+\Mm_X)$" to represent log discrepancies. This is because $(X,B+\Mm_X)$ is a sub-pair and the log discrepancies of this sub-pair may be different from the log discrepancies of the generalized pair $(X,B,\Mm)$.

\subsection{Equidimensional reduction}

\begin{thm}\label{thm: has19 weak semistable reduction}
Let $(X,B)$ be a dlt pair and $\pi: X\rightarrow U$ a projective surjective morphism over a normal variety $U$. Then there exists a commutative diagram of projective morphisms
\begin{center}$\xymatrix{
Y\ar@{->}[r]^{f}\ar@{->}[d]_{\pi'} & X\ar@{->}[d]^{\pi}\\
V\ar@{->}[r]^{\varphi} & U
}$
\end{center}
such that
\begin{enumerate}
    \item $f,\varphi$ are birational morphisms, $\pi'$ is an equidimensional contraction, $Y$ only has $\Qq$-factorial toroidal singularities, and $V$ is smooth, and
    \item there exist two $\Rr$-divisors $B_Y$ and $E$ on $Y$, such that
    \begin{enumerate}
    \item $K_Y+B_Y=f^*(K_X+B)+E$,
    \item $B_Y\geq 0$, $E\geq 0$, and $B_Y\wedge E=0$,
    \item $(Y,B_Y)$ is lc quasi-smooth, and any lc center of $(Y,B_Y)$ on $X$ is an lc center of $(X,B)$.
    \end{enumerate}
\end{enumerate}
\end{thm}
\begin{proof} This result follows from \cite{AK00}, see also \cite[Theorem B.6]{Hu20}, \cite[Theorem 2]{Kaw15} and \cite[Step 2 of the proof of Lemma 3.1]{Has19}.
\end{proof}

\subsection{Iitaka dimensions}

We refer the readers to \cite[Section 2]{HH20} for the formal definitions and basic properties of $\kappa_\sigma(X/U,D)$ and $\kappa_\iota(X/U,D)$.

\begin{lem}[{cf. \cite[V. 2.6(5) Remark]{Nak04}}]\label{lem: num dimension >=0 imply pe}
Let $X$ be a normal projective variety and $D$ an $\Rr$-Cartier $\Rr$-divisor on $X$ such that $\kappa_{\sigma}(D)\geq 0$. Then $D$ is pseudo-effective.
\end{lem}
\begin{proof}
By definition, there exists a Cartier divisor $A$ on $X$ such that $\sigma(D;A)\geq 0$. In particular, there exists a sequence of strictly increasing positive integers $m_i$, such that $\dim H^0(X,\lfloor m_iD\rfloor+A)>0$, hence $\lfloor m_iD\rfloor+A$ is effective for any $i$. Thus $m_iD+A$ is effective for any $i$, hence $D+\frac{1}{m_i}A$ is effective for any $i$. Thus $D$ is the limit of the effective $\Rr$-divisors $D+\frac{1}{m_i}A$, hence $D$ is pseudo-effective.
\end{proof}

\begin{lem}\label{lem: property of numerical and Iitaka dimension} Let $\pi: X\rightarrow U$ be a projective morphism from a normal variety to a variety, and $D$ an $\Rr$-Cartier $\Rr$-divisor on $X$. Then:
\begin{enumerate}
    \item $D$ is big$/U$ if and only if $\kappa_{\sigma}(X/U,D)=\dim X-\dim U$.
    \item Let $D_1,D_2$ be two $\Rr$-Cartier $\Rr$-divisors on $X$. Suppose that $D_1\sim_{\mathbb R,U}E_1\geq 0$ and $D_2\sim_{\mathbb R,U}E_2\geq 0$ for some $\Rr$-divisors $E_1,E_2$ such that $\Supp E_1=\Supp E_2$. Then $\kappa_{\sigma}(X/U,D_1)=\kappa_{\sigma}(X/U,D_2)$ and $\kappa_{\iota}(X/U,D_1)=\kappa_{\iota}(X/U,D_2)$.
    \item Let $f: Y\rightarrow X$ be a surjective birational morphism and $D_Y$ an $\Rr$-Cartier $\Rr$-divisor on $Y$ such that $D_Y=f^*D+E$ for some $f$-exceptional $\Rr$-divisor $E\geq 0$. Then $\kappa_{\sigma}(Y/U,D_Y)=\kappa_{\sigma}(X/U,D)$ and $\kappa_{\iota}(Y/U,D_Y)=\kappa_{\iota}(X/U,D)$. 
    \item Let $g: Z\rightarrow X$ be a surjective morphism from a normal variety such that $Z$ is projective over $U$. Then $\kappa_{\sigma}(Z/U,g^*D)=\kappa_{\sigma}(X/U,D)$ and $\kappa_{\iota}(Z/U,g^*D)=\kappa_{\iota}(X/U,D)$.
    \item  Let $\bar D$ be an $\Rr$-Cartier $\Rr$-divisor on $X$ such that $D\equiv_U\bar D$. Then $\kappa_{\sigma}(X/U,D)=\kappa_{\sigma}(X/U,\bar D)$.
    \item  Let $\phi: X\dashrightarrow X'$ be a partial $D$-MMP$/U$ and let $D':=\phi_*D$. Then $\kappa_{\sigma}(X/U,D)=\kappa_{\sigma}(X'/U,D')$ and $\kappa_{\iota}(X/U,D)=\kappa_{\iota}(X'/U,D')$
\end{enumerate}
\end{lem}
\begin{proof}
For (1)-(5), let $F$ be a very general fiber of the Stein factorization of $\pi$. Possibly replacing $X$ with $F$, $U$ with $\{pt\}$, and $D,D_1,D_2,\bar D$ with $D|_F,D_1|_F,D_2|_F,\bar D|_F$ respectively, we may assume that $X$ is projective and $U=\{pt\}$. (2) follows from \cite[Remark 2.8(1)]{HH20} and (3)(4) follow from \cite[Remark 2.8(2)]{HH20}.

To prove (1)(5), let $h: \tilde X\rightarrow X$ be a resolution of $X$. By (4), we may replace $X$ with $\tilde X$, $D$ with $h^*D$, and $\bar D$ with $h^*\bar D$, and assume that $X$ is smooth. 

If $D$ is big, then  $\kappa_{\sigma}(D)=\dim X$ by definition. If $\kappa_{\sigma}(D)=\dim X$, then $D$ is pseudo-effective by Lemma \ref{lem: num dimension >=0 imply pe}, hence $D$ is big by \cite[V. 2.7(3) Proposition]{Nak04}. This gives (1).

To prove (5), notice that $D$ is pseudo-effective if and only if $\bar D$ is pseudo-effective. If $D$ is not pseudo-effective, then $\kappa_{\sigma}(D)=\kappa_{\sigma}(\bar D)=-\infty$ by Lemma \ref{lem: num dimension >=0 imply pe}. If $D$ is pseudo-effective, then (5) follows from \cite[V. 2.7(1) Proposition]{Nak04}.

To prove (6), let $p: W\rightarrow X$ and $q: W\rightarrow X'$ be a common resolution such that $q=\phi\circ p$. Then $p^*D=q^*D'+F$ for some $F\geq 0$ that is $q$-exceptional. By (3), we have
$$\kappa_{\sigma}(X/U,D)=\kappa_{\sigma}(W/U,p^*D)=\kappa_{\sigma}(W/U,q^*D'+F)=\kappa_{\sigma}(X'/U,D')$$
and
\[\kappa_{\iota}(X/U,D)=\kappa_{\iota}(W/U,p^*D)=\kappa_{\iota}(W/U,q^*D'+F)=\kappa_{\iota}(X'/U,D').\qedhere\]
\end{proof}

\begin{lem}\label{lem: relative subaddivitiy iitaka dimensions}
Let $f: X\rightarrow Y$ and $g: Y\rightarrow Z$ be two contractions between normal quasi-projective varieties such that general fibers of $Y\rightarrow Z$ are smooth and $Y$ is $\Qq$-Gorenstein. Let $(X,B)$ be a pair that is lc over a non-empty open subset of $Y$. Let $D$ be an $\Rr$-Cartier $\Rr$-divisor on $X$ such that $D-(K_{X/Y}+B)$ is nef$/Z$. Then for any $\Rr$-Cartier $\Rr$-divisor $Q$ on $Y$, we have
$$\kappa_{\sigma}(X/Z,D+f^*Q)\geq\kappa_{\sigma}(X/Y,D)+\kappa(Y/Z,Q).$$
\end{lem}
\begin{proof}
Let $z\in Z$ be a very general point and let $X_z:=(g\circ f)^{-1}(z),Y_z:=g^{-1}(z)$ be the fibers of $X$ and $Y$ over $z$ respectively. We have an induced contraction $f_z: X_z\rightarrow Y_z$. Let $F$ be a very general fiber of $f_z$. Then $F$ is also a very general fiber of $f$. 

First assume that $\dim Y>\dim Z$. By our assumption, $Y_z$ is smooth, $(X_z,B|_{X_z})$ is lc over a non-empty open subset of $Y_z$, and
$$D|_{X_z}-(K_{X_z/Y_z}+B|_{X_z})=(D-(K_{X/Y}+B))|_{X_z}$$
is nef. By \cite[(3.3)]{Fuj20},
\begin{align*}
    \kappa_{\sigma}(X/Z,D+f^*Q)&=\kappa_{\sigma}(X_z,D|_{X_z}+f_z^*Q|_{Y_z})\geq \kappa_{\sigma}(X_z/Y_z,D|_{X_z})+\kappa(Y_z,Q|_{Y_z})\\
    &=\kappa_{\sigma}(F,D|_F)+\kappa(Y/Z,Q)=\kappa_{\sigma}(X/Y,D)+\kappa(Y/Z,Q).
\end{align*}

Now assume that $\dim Y=\dim Z$ so that $\kappa(Y/Z,Q)=0$. If $\dim X=\dim Y$ then there is nothing left to prove, so we may assume that $\dim X>\dim Y$. In this case $f^*Q|_{X_z}=0$, so we have
\begin{align*}
\kappa_{\sigma}(X/Z,D+f^*Q)&=\kappa_{\sigma}(X_z,D|_{X_z}+f^*Q|_{X_z})=\kappa_{\sigma}(X_z,D|_{X_z})= \kappa_{\sigma}(X/Z,D)\\
    &=\kappa_{\sigma}(X/Y,D)=\kappa_{\sigma}(X/Y,D)+\kappa(Y/Z,Q),
\end{align*}
and we are done.
\end{proof}

\begin{lem}\label{lem: iitaka fibration numerical abundant divisor gpair dimension}
Let $(X,B,\Mm)/U$ be an lc g-pair such that $K_X+B+\Mm_X\equiv_{U}G$ for some $\Rr$-divisor $G\geq 0$, such that $U$ is quasi-projective and $G$ is abundant over $U$. Let $X\dashrightarrow V$ be the Iitaka fibration over $U$ associated to $G$, and $(W,B_W,\Mm)$ a log smooth model of $(X,B,\Mm)$ such that the induced map $\psi: W\rightarrow V$ is a morphism over $U$. Then
\begin{enumerate}
    \item $\kappa_{\sigma}(W/U,K_W+B_W+\Mm_W)=\dim V-\dim U$, and
    \item $\kappa_{\sigma}(W/V,K_W+B_W+\Mm_W)=0$.
\end{enumerate}
\end{lem}
\begin{proof}
Let $h_V: \bar V\rightarrow V$ be a resolution of $V$. By Lemmas \ref{lem: property of numerical and Iitaka dimension}(3) and \cite[Lemma 3.6]{HL21a} possibly replacing $(W,B_W,\Mm)/U$ with a higher model, we may assume that the induced map $\bar\psi: W\rightarrow\bar V$ is a morphism. Since $(W,B_W,\Mm)$ a log smooth model of $(X,B,\Mm)$, we have
$$K_W+B_W+\Mm_W=h^*(K_X+B+\Mm_X)+E$$ where $h: W\to X$ is the induced morphism, $\Mm$ descends to $W$, and $E\geq 0$.
\begin{center}$\xymatrix{
W\ar@{->}[d]_{\bar\psi}\ar@{->}[drr]^{\psi}\ar@{->}[rr]^{h} &  & X\ar@{-->}[d] \\
\bar V\ar@{->}[dr]\ar@{->}[rr]^{h_V}&    & V\ar@{->}[dl] \\
 & U &
}$
\end{center}
Since $G\geq 0$ is abundant over $U$, by \cite[Proposition 2.2.2(1)]{Cho08}, $$\dim V-\dim U=\kappa(X/U,G)=\kappa_{\iota}(X/U,G)=\kappa_{\sigma}(X/U,G)\geq 0.$$
Since $X\dashrightarrow V$ is the Iitaka fibration associated to $G$ over $U$, there exists an effective ample$/U$ $\Rr$-divisor $A$ on $V$ and an $\Rr$-divisor $F\geq 0$ on $W$ such that $h^*G=\psi^*A+F$ 
for some $h$-exceptional $\Rr$-divisor $F\geq 0$ on $W$. Then for any real number $k$, we have
$$K_W+B_W+\Mm_W+k\psi^*A\equiv_{U}(1+k)\psi^*A+E+F.$$
By Lemma \ref{lem: property of numerical and Iitaka dimension}(2)(3)(5), for any $k\geq 0$ we have \begin{align*}
    \kappa_{\sigma}(W/U,K_W+B_W+\Mm_W+k\psi^*A)&=\kappa_{\sigma}(W/U,(1+k)\psi^*A+E+F)=\kappa_{\sigma}(W/U,\psi^*A+E+F)\\
    &=\kappa_{\sigma}(W/U,K_W+B_W+\Mm_W)=\kappa_{\sigma}(X/U,K_X+B+\Mm_X)\\
    &=\kappa_{\sigma}(X/U,G)=\kappa(X/U,G)=\dim V-\dim U.
\end{align*}
In particular, we get (1). Since $A$ is ample$/U$, $h_V^*A$ is big$/U$, and we may pick a sufficiently large positive integer $k$ such that $K_{\bar V}+kh_V^*A$ is big$/U$.

Since $(W,B_W,\Mm)$ is a log smooth model of $(X,B,\Mm)$, $(W,B_W)$ is lc. Since $\bar V$ is smooth, any very general fiber of the induced morphism $\bar V\rightarrow U$ is smooth. Let $D:=K_W+B_W+\Mm_W-\bar\psi^*K_{\bar V}$ and $Q:=K_{\bar V}+kh_V^*A$. Then $D-(K_{W/\bar V}+B_W)=\Mm_W$ is nef$/U$. By Lemma \ref{lem: relative subaddivitiy iitaka dimensions} and noticing that the restriction of $\bar\psi^*K_{\bar V}$ to a general fiber of $\bar\psi$ is zero, we have
 \begin{align*}
    \dim V-\dim U&=\kappa_{\sigma}(W/U,K_W+B_W+\Mm_W+k\psi^*A)=\kappa_{\sigma}(W/U,K_W+B_W+\Mm_W+k\bar\psi^*h_V^*A)\\
    &=\kappa_{\sigma}(W/U,D+\bar\psi^*Q)\geq \kappa_{\sigma}(W/\bar V,D)+\kappa(\bar V/U,Q)\\
    &=\kappa_{\sigma}(W/\bar V,K_W+B_W+\Mm_W-\bar\psi^*K_{\bar V})+\kappa(\bar V/U,K_{\bar V}+kh_V^*A)\\
    &=\kappa_{\sigma}(W/\bar V,K_W+B_W+\Mm_W)+(\dim V-\dim U).
\end{align*}
Thus $\kappa_{\sigma}(W/\bar V,K_W+B_W+\Mm_W)\leq 0$, hence $\kappa_{\sigma}(W/V,K_W+B_W+\Mm_W)\leq 0$. Since $K_W+B_W+\Mm_W\equiv_U h^*G+E\geq 0$, $\kappa_{\sigma}(W/V,K_W+B_W+\Mm_W)\geq 0$. Thus $\kappa_{\sigma}(W/V,K_W+B_W+\Mm_W)=0$, and we get (2).
\end{proof}

\subsection{Preliminaries on the MMP for generalized pairs}

\begin{lem}[{\cite[Lemma 2.20]{HL21a}, cf. \cite[Proposition 3.9]{HL22}}]\label{lem: rlinear version of HL22 3.8}
Let $(X,B,\Mm)/U$ be a $\Qq$-factorial lc g-pair such that $X$ is klt and $K_X+B+\Mm_X\equiv_{U}D_1-D_2$ (resp.  $\sim_{\mathbb R,U}D_1-D_2$) where $D_1\geq 0$, $D_2\geq 0$ have no common components. Suppose that $D_1$ is very exceptional over $U$ (see \cite[Definition 3.1]{Bir12}). Then any $(K_X+B+\Mm_X)$-MMP$/U$ with scaling of an ample$/U$ $\Rr$-divisor either terminates with a Mori fiber space or contracts $D_1$ after finitely many steps. Moreover, if $D_2=0$, then this MMP terminates with a model $Y$ such that  $K_Y+B_Y+\Mm_Y\equiv_{U}0$ (resp. $\sim_{\mathbb R,U}0$), where $B_Y$ is the strict transform of $B$ on $Y$. 
\end{lem}


\begin{lem}[{\cite[Lemma 2.25]{HL21a}}]\label{lem: still an mmp under perturbation}
Let $X\rightarrow U$ be a projective morphism such that $X$ is normal quasi-projective. Let $D,A$ be two $\Rr$-Cartier $\Rr$-divisors on $X$ and let $\phi: X\dashrightarrow X'$ be a partial $D$-MMP$/U$. Then there exists a positive real number $t_0$, such that for any $t\in (0,t_0]$, $\phi$ is also a partial $(D+tA)$-MMP$/U$. Note that $A$ is not necessarily effective.
\end{lem}
\begin{proof}
We let
$$X:=X_0\dashrightarrow X_1\dashrightarrow\dots\dashrightarrow X_n=X'$$
be this partial MMP, and $D_i,A_i$ the strict transforms of $D$ and $A$ on $X_i$ respectively. Let $X_i\rightarrow Z_i$ be the $D_i$-negative extremal contraction of a $D_i$-negative extremal ray $R_i$ in this MMP for each $i$. Then $D_i\cdot R_i<0$ for each $i$. Thus there exists a positive real number $t_0$ such that $(D_i+t_0A_i)\cdot R_i<0$ for each $i$. In particular, $(D_i+tA_i)\cdot R_i<0$ for any $i$ and any $t\in (0,t_0]$. Thus $\phi$ is a partial $(D+tA)$-MMP$/U$ for any $t\in (0,t_0]$.
\end{proof}

\begin{lem}[{cf. \cite[Lemma 2.17]{LT22a}}]\label{lem: limit movable r divisors gpairs}
Let $(X,B,\Mm)/U$ be a $\Qq$-factorial NQC lc g-pair such that $X$ is klt and $K_X+B+\Mm_X$ is pseudo-effective$/U$. Let $A\geq 0$ be an ample$/U$ $\Rr$-divisor on $X$ such that $(X,B+A,\Mm)$ is lc and $K_X+B+A+\Mm_X$ is nef$/U$. Let $$(X,B,\Mm):=(X_0,B_0,\Mm)\dashrightarrow (X_1,B_1,\Mm)\dashrightarrow\dots\dashrightarrow (X_i,B_i,\Mm)\dashrightarrow\dots$$
be a $(K_X+B+\Mm_X)$-MMP$/U$ with scaling of $A$, and $A_i$ the strict transform of $A$ on $X_i$ for each $i$. Then there exists a positive integer $n$ and a positive real number $\epsilon_0$, such that $K_{X_j}+B_j+\epsilon A_j+\Mm_{X_j}$ is movable$/U$ for any $\epsilon\in (0,\epsilon_0)$ and $j\geq n$. In particular, $K_{X_j}+B_j+\Mm_{X_j}$ is a movable$/U$ (cf. Definition \ref{defn: rel nz decomposition}) $\Rr$-divisor  for any $j\ge n$.
\end{lem}
\begin{proof}
Let $\lambda_i$ be the $i$-th scaling number of this MMP for each $i$, i.e.
$$\lambda_i:=\inf\{t\geq 0 \mid K_{X_i}+B_i+tA_i+\Mm_{X_i}\text{ is nef/}U\}.$$
We may assume that this MMP does not terminate. By \cite[Theorem 2.24]{HL21a}, we have $\lim_{i\rightarrow+\infty}\lambda_i=0$.

Let $n$ be the minimal positive integer such that $X_{i}\dashrightarrow X_{i+1}$ is a flip for any $i\geq n$. If $\lambda_i<\lambda_{i-1}$, then $X\dashrightarrow X_{i}$ is a $(K_X+B+tA+\Mm_X)$-MMP$/U$ with scaling of $(1-t)A$ for any $t\in [\lambda_i,\lambda_{i-1})$. Since $X$ is $\Qq$-factorial klt, there exists $\Delta_t\sim_{\Rr,U}B+tA+\Mm_X$ such that $(X,\Delta_t)$ is klt and $\Delta_t$ is big for any $t\in (0,1]$. By \cite[Corollary 3.9.2]{BCHM10}, $K_{X_i}+B_i+tA_i+\Mm_{X_i}$ is semi-ample$/U$ for any $i$ and any $t\in [\lambda_i,\lambda_{i-1})$. Let $\epsilon_0:=\lambda_n$. Then for any $\epsilon\in (0,\epsilon_0)$, there exists $i\geq n$ such that $\lambda_i<\lambda_{i-1}$ and $\epsilon\in [\lambda_i,\lambda_{i-1})$, and $K_{X_i}+B_i+\epsilon A_i+\Mm_{X_i}$ is semi-ample$/U$. Since $X_i\dashrightarrow X_j$ is small for any $i,j\geq n$, $K_{X_j}+B_j+\epsilon A_j+\Mm_{X_j}$ is movable$/U$ for any $j\geq n$ and $\epsilon\in (0,\epsilon_0)$, and $K_{X_j}+B_j+\Mm_{X_j}$ is a $/U$ $\Rr$-divisor.
\end{proof}

\begin{lem}\label{lem: limit of movable divisors mmp only contain flips}
Let $X\rightarrow U$ be a projective morphism such that $X$ is quasi-projective. Assume that $D$ is an $\Rr$-Cartier $\Rr$-divisor on $X$ such that $D$ is a movable$/U$ $\Rr$-divisor on $X$, and let $\phi: X\dashrightarrow X'$ be a partial $D$-MMP$/U$. Then $\phi$ only contains flips.
\end{lem}
\begin{proof}
Since $D$ is a movable$/U$ $\Rr$-divisor, $D$ is pseudo-effective$/U$, so $\phi$ only contains flips and divisorial contractions. 

If $\phi$ contains a divisorial contraction, let $\psi: X_1\rightarrow X_1'$ be the first divisorial contraction in $\phi$. Let $D_1$ be the strict transform of $D$ on $X_1$. Then $X\dashrightarrow X_1$ only contains flips, hence it is an isomorphism in codimension one, so $D_1$ is also a movable$/U$ $\Rr$-divisor on $X_1$. Let $D_1':=\psi_*D_1$. Then
$$D_1=\psi^*D_1'+F$$
for some $F\geq 0$ that is exceptional over $X_1'$.

Since $D_1$ is a movable$/U$ divisor,  $D_1'$ is a movable$/X_1'$ divisor. 
Thus for any very general $\psi$-exceptional curve $C$, $D_1\cdot C\geq 0$. By the general negativity lemma \cite[Lemma 3.3]{Bir12}, $-F\geq 0$. Thus $F=0$, and $\psi$ cannot be a $D_1$-negative extremal contraction, a contradiction. Thus $\phi$ only contains flips.
\end{proof}

\begin{lem}\label{lem: gmmp scaling numbers go to 0}
Let $(X,B,\Mm)/U$ be a $\Qq$-factorial NQC lc g-pair. Let $H\geq 0$ be an $\Rr$-divisor on $X$ such that $(X,B+H,\Mm)$ is lc and $K_X+B+H+\Mm_X$ is nef$/U$. Assume that $(X,B+\mu H,\Mm)/U$ has a log minimal model for any $\mu\in (0,1]$. Then we can construct a $(K_X+B+\Mm_X)$-MMP$/U$ with scaling of $H$:
$$(X,B,\Mm):=(X_0,B_0,\Mm)\dashrightarrow (X_1,B_1,\Mm)\dashrightarrow\dots\dashrightarrow (X_i,B_i,\Mm)\dashrightarrow\dots.$$
Let $H_i$ be the strict transform of $H$ on $X_i$ for each $i$, and let
$$\lambda_i:=\inf\{t\geq 0 \mid K_{X_i}+B_i+tH_i+\Mm_{X_i}\text{ is nef/}U\}$$
be the $i$-th scaling number of this MMP for each $i$. Then this MMP
\begin{enumerate}
    \item either terminates after finitely many steps, or
    \item does not terminate and $\lim_{i\rightarrow+\infty}\lambda_i=0$.
\end{enumerate}
\end{lem}
\begin{proof}
If $\lambda_0=0$ then there is nothing left to prove. So we may assume that $\lambda_0>0$. By \cite[Lemma 3.21]{HL22}, we may pick $\lambda_0'\in (0,\lambda_0)$ such that any sequence of the $(K_X+B+\lambda_0'H+\Mm_X)$-MMP$/U$ is $(K_X+B+\lambda_0H+\Mm_X)$-trivial. 

By \cite[Theorem 2.24]{HL21a}, we may run a $(K_X+B+\lambda_0'H+\Mm_X)$-MMP$/U$ with scaling of a general ample$/U$ divisor, which terminates with a log minimal model. We let
$$(X,B,\Mm):=(X_0,B_0,\Mm)\dashrightarrow (X_1,B_1,\Mm)\dashrightarrow\dots\dashrightarrow (X_{k_1},B_{k_1},\Mm)$$
be this sequence of the MMP$/U$. Then this sequence consists of finitely many steps of a $(K_X+B+\Mm_X)$-MMP$/U$ with scaling of $H$, with scaling numbers $\lambda_0=\lambda_1=\dots=\lambda_{k_1-1}$. Since
$$K_{X_{k_1}}+B_{k_1}+\lambda'_1H_{k_1}+\Mm_{X_{k_1}}$$
is nef$/U$, we have $\lambda_{k_1}\leq\lambda_1'<\lambda_1$. 

We may replace $(X,B,\Mm)/U$ with $(X_{k_1},B_{k_1},\Mm)/U$ and continue this process. If this MMP does not terminate, then we may let $\lambda:=\lim_{i\rightarrow+\infty}\lambda_i$. By our construction, $\lambda\not=\lambda_i$ for any $i$, and the lemma follows from \cite[Theorem 4.1]{LT22b}. 
\end{proof}

\begin{lem}[{cf. \cite[3.5 Lifting flips, Page 727-728]{HL22}, \cite[2.5 Lifting a sequence of flips with scaling, Lemma 2.13]{LT22b}}]\label{lem: lift mmp}
Let $(X,B,\Mm)/U$ be an NQC lc g-pair, $S$ an lc center of $(X,B,\Mm)$,  $(Y,B_Y,\Mm)$ a dlt model of $(X,B,\Mm)$ with induced birational morphism $f: Y\rightarrow X$, and $S_Y$ a component of $\lfloor B_Y\rfloor$ such that $f(S_Y)=S$. Let
$$\phi: (X,B,\Mm)\dashrightarrow (X',B',\Mm)$$
be a partial $(K_X+B+\Mm_X)$-MMP$/U$ and $S'$ an lc center of $(X',B',\Mm)$ such that $\phi|_S: S\dashrightarrow S'$ is a birational map. Then there exists a partial $(K_Y+B_Y+\Mm_Y)$-MMP$/U$
$$\psi: (Y,B_Y,\Mm)\dashrightarrow (Y',B_Y',\Mm),$$
such that
\begin{enumerate}
    \item $(Y',B_Y',\Mm)$ is a dlt model of $(X',B',\Mm)$, and
    \item the strict transform of $S_Y$ on $Y'$ is a component of $\lfloor B_Y'\rfloor$. 
\end{enumerate}
\end{lem}
\begin{proof}
We only need to prove the lemma when $\phi$ is a divisorial contraction or a flip. If $\phi$ is a flip, then we let $X\rightarrow Z$ be the flipping contraction and let $X'\rightarrow Z$ be the flipped contraction.  The rest of the proof of (1) is similar to the \cite[First paragraph of the proof of Lemma 2.13]{LT22b}: If $\phi$ is a divisorial contraction, then we let $Z=X'$. Thus $(X',B',\Mm)/Z$ is a log minimal model of $(X,B,\Mm)/Z$ such that $K_{X'}+B'+\Mm_{X'}$ is ample$/Z$. By \cite[Lemmas 3.9, 3.15]{HL21a} and \cite[Theorem 3.14]{HL21a}, we may run a $(K_{Y}+B_Y+\Mm_Y)$-MMP$/Z$ with scaling of an ample$/Z$ divisor which terminates with a good minimal model $(Y',B_Y',\Mm)/Z$. By \cite[Lemma 3.9]{HL21a}, $(Y',B_Y',\Mm)$ is a dlt model of $(X',B',\Mm)$, and we get (1).

We let $p: W\rightarrow Y$ and $q: W\rightarrow Y'$ be a resolution of indeterminacies of the induced birational map $\phi_Y: Y\dashrightarrow Y'$. By \cite[Lemma 3.8]{HL21a}, $p^*(K_Y+B_Y+\Mm_Y)=q^*(K_{Y'}+B_{Y'}+\Mm_{Y'})+F$ where $F\geq 0$ is exceptional$/Y'$, and $\Supp p_*F$ contains all $\phi_Y$-exceptional divisors. By (1), $a(S_Y,Y,B_{Y'},\Mm)=0$, hence $S_Y$ is not a component of $\Supp p_*F$, and we get (2).
\end{proof}

\subsection{Proper log smooth models}

\begin{defn}[Log smooth model, {\cite[Definition 3.1]{HL21a}}]\label{defn: log smooth models}
Let $(X,B,\Mm)/U$ be an lc g-pair and $h: W\rightarrow X$ a log resolution of $(X,\Supp B)$ such that $\Mm$ descends to $W$. Let $B_W\geq 0$ and $E\geq 0$ be two $\Rr$-divisors on $W$ such that
\begin{enumerate}
    \item $K_W+B_W+\Mm_W=h^*(K_X+B+\Mm_X)+E$,
    \item $(W,B_W)$ is log smooth dlt,
    \item $E$ is $h$-exceptional, and
    \item for any $h$-exceptional prime divisor $D$ such that $a(D,X,B,\Mm)>0$, $D$ is a component of $E$.
\end{enumerate}
Then $(W,B_W,\Mm)$ is called a \emph{log smooth model} of $(X,B,\Mm)$. If we additionally assume that
\begin{itemize}
    \item[(5)] for any $h$-exceptional prime divisor $D$ such that $a(D,X,B,\Mm)>0$, $D$ is a component of $\{B_W\}$,
\end{itemize}
then $(W,B_W,\Mm)$ is called a \emph{proper log smooth model} of $(X,B,\Mm)$.
\end{defn}

\begin{lem}\label{lem: special proper log smooth model}
Let $(X,B,\Mm)/U$ be an lc g-pair. Then there exists a proper log smooth model $(W,B_W=B_W^h+B_W^v,\Mm)$ of $(X,B,\Mm)$, such that
\begin{enumerate}
    \item $B_W^h\geq 0$ and $B_W^v$ is reduced,
    \item $B_W^v$ is vertical over $U$, and
    \item for any real number $t\in (0,1]$, all lc centers of $(W,B_W-tB_W^v,\Mm)$ dominate $U$.
\end{enumerate}

\end{lem}
\begin{proof}
By \cite[Lemma 3.6]{HL21a}, possibly replacing $(X,B,\Mm)$ with a proper log smooth model, we may assume that $(X,\Supp B)$ is log smooth and $\Mm$ descends to $X$. By \cite[Lemma 2.10]{Has18}, there exists a proper log smooth model $(W,B_W=B_W^h+B_W^v)$ of $(X,B)$, such that
\begin{itemize}
    \item $B_W^h\geq 0$ and $B_W^v$ is reduced,
    \item $B_W^v$ is vertical over $U$, and
    \item for any real number $t\in (0,1]$, all lc centers of $(W,B_W-tB_W^v)$ dominate $U$.
\end{itemize}
Since $\Mm$ descends to $X$, $(W,B_W,\Mm)$ is a proper log smooth model of $(X,B,\Mm)$, and for any real number $t\in (0,1]$, any lc center of $(W,B_W-tB_W^v,\Mm)$ is an lc center of $(W,B_W-tB_W^v)$ and dominates $U$. Thus $(W,B_W=B_W^h+B_W^v,\Mm)$ satisfies our requirements.
\end{proof}

\subsection{Canonical bundle formula}

\begin{thm}\label{thm: numerical generalized canonical bundle formula}
Let $(X,B,\Mm)/U$ be an NQC lc g-pair such that $U$ is quasi-projective, and let $\pi: X\rightarrow V$ be a surjective morphism over $U$. Assume that $K_X+B+\Mm_X\sim_{\mathbb R,V}0$. Then there exists an NQC lc g-pair $(V,B_V,\Mm^V)/U$, such that
\begin{enumerate}
    \item  $K_X+B+\Mm_X\sim_{\mathbb R}\pi^*(K_V+B_V+\Mm^V_V)$,
    \item any lc center of $(V,B_V,\Mm^V)$ is the image of an lc center of $(X,B,\Mm)$ in $V$, and
    \item if all lc centers of $(X,B,\Mm)$ dominate $V$, then $(V,B_V,\Mm^V)$ is klt. 
\end{enumerate}
\end{thm}
\begin{proof}
By the theory of Shokurov-type rational polytopes (cf. \cite[Proposition 3.20]{HL22}) and the theory of uniform rational polytopes (see \cite[Lemma 5.3]{HLS19}, \cite[Therem 1.4]{Che20}), we may assume that $(X,B,\Mm)/U$ is a $\Qq$-g-pair. 

\medskip

\noindent\textbf{Step 1}. In this step we prove the case when $X\rightarrow V$ is a generically finite morphism. 

By \cite[Theorem 4.5, (4.3),(4.4)]{HL20}, there exists an lc $\Qq$-g-pair $(V,B_V,\Mm^V)/U$, such that $K_X+B+\Mm_X\sim_{\mathbb Q}\pi^*(K_V+B_V+\Mm^V_V)$, and $B_V$ and $\Mm^V$ are defined in the following way: 

Let $V^0$ be the smooth locus of $V$, $X^0:=X\times_VV^0$, and $\pi|_{X^0}: X^0\rightarrow V^0$ the restriction of $\pi$. Then we have the Hurwitz formula
$$K_{X^0}=(\pi|_{X^0})^*K_{V^0}+R^0$$
where $R^0$ is the effective ramification divisor of $f|_{X^0}$. Let $R$ be the closure of $R^0$ in $X$, and let $B_V:=\frac{1}{\deg \pi}\pi_*(R+B)$. For any proper birational morphism $\mu: V'\rightarrow V$, let $X'$ be the main component of $X\times_VV'$ with induced birational map $\pi': X'\rightarrow V'$. We let $\Mm^V_{V'}=\frac{1}{\deg\pi}\pi'_*\Mm_{X'}$.

(1) follows immediately. 

Since $(V,B,\Mm^V)/U$ is a g-pair, for any prime divisor $E$ over $V$, there exists a birational morphism $h_V: \tilde V\rightarrow V$ such that $\Mm^V$ descends to $\tilde V$ and $E$ is on $\tilde V$. We let $h: \tilde X\rightarrow X$ be a birational morphism such that $\Mm$ descends to $\tilde X$ and the induced map $\tilde\pi: \tilde X\rightarrow\tilde V$ is a  morphism. 
\begin{center}$\xymatrix{
X'\ar@{->}[d]_{\pi'}\ar@{->}[r]& X\ar@{->}[d]_{\pi} & \tilde X\ar@{->}[l]_{h}\ar@{->}[d]^{\tilde\pi}\\
V'\ar@{->}[r]^{\mu}& V & \tilde V\ar@{->}[l]_{h_V}
}$
\end{center}
There are two cases:

\medskip

\noindent\textbf{Case 1}. $E$ is exceptional over $V$. In this case we let $F\subset\tilde\pi^{-1}(E)$ be a prime divisor, and let $r\leq\deg f$ be the ramification index of $\tilde\pi$ along $F$. Near the generic point of $F$, we have
\begin{align*}
K_{\tilde X}&=~h^*(K_X+B+\Mm_X)+(a(F,X,B,\Mm)-1)F\\
&\sim_{\Qq}h^*\pi^*(K_V+B_V+\Mm_V)+(a(F,X,B,\Mm)-1)F    
\end{align*}

and
\begin{align*}
 K_{\tilde X}&=\tilde\pi^*K_{\tilde V}+(r-1)F\\
 &=\tilde\pi^*h_V^*(K_V+B_V+\Mm_V)+r(a(E,V,B_V,\Mm^V)-1)F+(r-1)F\\
 &=h^*\pi^*(K_V+B_V+\Mm_V)+(ra(E,V,B_V,\Mm^V)-1)F.
\end{align*}
Let $\tilde X\rightarrow \bar X\rightarrow V$ be the Stein factorization of $\pi\circ h=h_V\circ\tilde\pi$. Since $E$ is exceptional over $V$, $F$ is exceptional over $\bar X$. Therefore $aF\sim_{\Qq,\bar X}0$ iff $a=0$ (applying negativity lemma to both $aF$ and $-aF$).  By comparing the two expressions of $K_{\tilde X}$ above, we have
$$a(F,X,B,\Mm)-1=ra(E,V,B_V,\Mm^V)-1,$$
hence $a(F,X,B,\Mm)\geq 0$ if and only if $a(E,V,B_V,\Mm^V)\geq 0$ and $a(F,X,B,\Mm)>0$ if and only if $a(E,V,B_V,\Mm^V)>0$. Moreover, since $F\subset\tilde\pi^{-1}(E)$, if $E$ is an lc place of $(V,B_V,\Mm^V)$, then $F$ is an lc place of $(X,B,\Mm)$ and $\Center_VE$ is contained in the image of $\Center_XF$ in $V$. 

\medskip

\noindent\textbf{Case 2}. $E$ is not exceptional over $V$. In this case, if $E$ is not a component of $B_V$, then $a(E,V,B_V,\Mm^V)=1>0$. If $E$ is a component of $B_V$, then we may let $B_1,\dots,B_m\subset\pi^{-1}(E)$ be the prime divisors on $X$ lying over $V$ and let $d_i$ be the degree of the induced morphism $\pi|_{B_i}: B_i\rightarrow E$. By our construction of $B_V$, 
$$a(E,V,B_V,\Mm^V)=1-\mult_EB_V=1-\frac{\sum_{i=1}^md_i\mult_{B_i}B}{\deg\pi}.$$
Since $\sum_{i=1}^m d_i\leq\deg\pi$, $a(E,V,B_V,\Mm^V)\geq 0$ if $\mult_{B_i}B\leq 1$ for each $i$, and $a(E,V,B_V,\Mm^V)>0$ if $\mult_{B_i}B<1$ for each $i$. Moreover, since $B_i\subset\pi^{-1}(E)$ for each $i$, if $E$ is an lc place of $(V,B_V,\Mm^V)$, then $B_i$ is an lc place of $(X,B,\Mm)$ for some $i$ and $E$ is contained in the image of $B_i$ in $V$.

\medskip

By our discussions above, we finish the proof in the case when $X\rightarrow V$ is a generically finite morphism.

\medskip

\noindent\textbf{Step 2}. In this step we prove the case when $X\rightarrow V$ is a contraction. 

By \cite[Theorem 2.20]{FS20}, there exists an lc $\Qq$-g-pair $(V,B_V,\Mm^V)/U$, such that $K_X+B+\Mm_X\sim_{\mathbb Q}\pi^*(K_V+B_V+\Mm^V_V)$. Moreover, for any birational morphism $h_V: \tilde V\rightarrow V$, we have an $\Rr$-divisor $B_{\tilde V}$ satisfying $K_{\tilde V}+B_{\tilde V}+\Mm^V_{\tilde V}=h_V^*(K_V+B_V+\Mm^V_V)$ and defined in the following way: let $\tilde X$ be the main component of $X\times_{V}\tilde V$, and $h: \tilde X\rightarrow X$ and $\tilde\pi: \tilde X\rightarrow\tilde V$ the induced morphisms. Let $K_{\tilde X}+\tilde B+\Mm_{\tilde X}:=h^*(K_X+B+\Mm_X)$. For any prime divisor $E$ on $\tilde V$, $\mult_{E}B_{\tilde V}=1-t_E$, where $$t_E:=\sup\{s\mid (\tilde X,\tilde B+s\tilde \pi^*E,\Mm)\text{ is lc over the generic point of }E\}.$$
Note that $E$ may not be $\Qq$-Cartier but $\tilde \pi^*E$ is always defined over the generic point of $E$. 

(1) follows immediately. 

If $E$ is an lc place of $(V,B_V,\Mm^V)$ on $\tilde V$, then $t_E=0$, hence $\tilde \pi^*E$ contains an lc center $F$ of $(\tilde X,\tilde B,\Mm)$ over the generic point of $E$. We have $F\subset\Supp\tilde\pi^*E$ and $\tilde\pi(F)\subset E$, hence $\tilde\pi(F)=E$. Thus $E$ is the image of an lc center of $(\tilde X,\tilde B,\Mm)$ on $\tilde V$, hence $\Center_{V}E$ is the image of an lc center of $(X,B,\Mm)$ in $V$. 

By our discussions above, we finish the proof in the case when $X\rightarrow V$ is a contraction.

\medskip

\noindent\textbf{Step 3}. In this step we prove the general case. 

We let $X\xrightarrow{f}Y\xrightarrow{g}V$ be the Stein factorization of $\pi$. Then $K_X+B+\Mm_X\sim_{\Qq,Y}0$, $f: X\rightarrow Y$ is a contraction and $g: Y\rightarrow V$ is a finite morphism. By Step 2, $K_X+B+\Mm_X\sim_{\Qq}f^*(K_Y+B_Y+\Mm^Y_Y)$ for some lc $\Qq$-g-pair $(Y,B_Y,\Mm^Y)/U$ such that any lc center of $(Y,B_Y,\Mm^Y)$ is the image of an lc center of $(X,B,\Mm)$ in $Y$. Moreover, $K_Y+B_Y+\Mm^Y_Y\sim_{\Qq,V}0$. By Step 1, $K_Y+B_Y+\Mm^Y_Y\sim_{\Qq}g^*(K_V+B_V+\Mm^V_V)$ for some lc g-pair $(V,B_V,\Mm^V)/U$ such that any lc center of $(V,B_V,\Mm^V)$ is the image of an lc center of $(Y,B_Y,\Mm^Y)$ in $V$, hence the image of an lc center of $(X,B,\Mm)$ in $V$. We immediately get (1)(2) and (3) follows from (2).
\end{proof}

\subsection{Special termination}

\begin{defn}\label{defn: gpair set for difficulty}
Let $\Ii\subset [0,1]$ and $\Ii'\subset [0,+\infty)$ be two sets. We define
$$\mathbb S(\Ii,\Ii'):=\{1-\frac{1}{m}+\sum_j\frac{r_jb_j}{m}+\sum_i\frac{s_i\mu_i}{m}\mid m\in\mathbb N^+,r_i,s_i\in\mathbb N, b_j\in\Ii,\mu_j\in\Ii'\}\cap ( 0,1].$$
\end{defn}

\begin{prop}[{\cite[Proposition 2.10]{HL22}}]\label{prop: HL22 2.8}
Let $\Ii\subset [0,1]$ and $\Ii'\subset [0,+\infty)$ be two sets. Let $(X,B,\Mm)/U$ be a $\Qq$-factorial NQC dlt g-pair such that $B\in\Ii$ and $\Mm=\sum\mu_i\Mm_i$, where $\mu_i\in\Ii'$ for each $i$ and each $\Mm_i$ is nef$/U$ $\bb$-Cartier. Then for any lc center $S$ of $(X,B,\Mm)$, the g-pair $(S,B_S,\Mm^S)/U$ given by the adjunction
$$K_S+B_S+\Mm^S_S:=(K_X+B+\Mm_X)|_S$$
is dlt, and $B_S\in\mathbb S(\Ii,\Ii')$.
\end{prop}

\begin{defn}[Difficulty, {\cite[Definition 4.5]{HL22}}]\label{defn: gpair difficulty}
Let $\Ii$ and $\Ii'$ be two finite sets of non-negative real numbers. Let $(X,B,\Mm)/U$ be a $\Qq$-factorial NQC dlt g-pair such that $B\in\Ii$ and $\Mm=\sum\mu_i\Mm_i$, where $\mu_i\in\Ii'$ for each $i$ and each $\Mm_i$ is nef$/U$ $\bb$-Cartier. For any lc center $S$ of $(X,B,\Mm)$ of dimension $\geq 1$, let $(S,B_S,\Mm^S)$ be the g-pair given by the generalized adjunction
$$K_S+B_S+\Mm^S_S:=(K_X+B+\Mm_X)|_S,$$
then we define
\begin{align*}
    d_{\Ii,\Ii'}(S,B_S,\Mm^S):=&\sum_{\alpha\in\mathbb S(\Ii,\Ii')}\#\{E\mid a(E,B_S,\Mm^S)<1-\alpha,\Center_SE\not\subset\lfloor B_S\rfloor\}\\
    &+\sum_{\alpha\in\mathbb S(\Ii,\Ii')}\#\{E\mid a(E,B_S,\Mm^S)\leq 1-\alpha,\Center_SE\not\subset\lfloor B_S\rfloor\}.
\end{align*}
\end{defn}

The following special termination result is similar to \cite{Fuj07,LMT20,HL22}. The proofs are also similar. For the reader's convenience, we provide a full proof here.

\begin{lem}\label{lem: special termination reduce to flip lemma}
Let $(X,B,\Mm)/U$ be a $\Qq$-factorial NQC dlt g-pair and let
$$(X,B,\Mm):=(X_0,B_0,\Mm)\dashrightarrow (X_1,B_1,\Mm)\dashrightarrow\dots\dashrightarrow (X_i,B_i,\Mm)\dashrightarrow\dots$$
be a $(K_X+B+\Mm_X)$-MMP$/U$. Let $\phi_{i,j}: X_i\dashrightarrow X_{j}$ be the induced birational maps for each $i$. For any $i\geq 0$ and any lc center $S_i$ of $(X_i,B_i,\Mm)$ of dimension $\geq 1$, we let $(S_i,B_{S_i},\Mm^{S_i})/U$ be the generalized pair given by adjunction
$$K_{S_i}+B_{S_i}+\Mm^{S_i}_{S_i}:=(K_{X_i}+B_i+\Mm_{X_i})|_{S_i}.$$
Then we have the following.
\begin{enumerate}
\item For any $i\gg 0$, $j\geq i$, and any lc center $S_i$ of $(X_i,B_i,\Mm)$, $\phi_{i,j}$ induces an isomorphism near the generic point of $S_i$. In particular, for any $i,j\gg 0$ and any lc center $S_i$ of $(X_i,B_i,\Mm)$, we may let $S_{i,j}$ be the strict transform of $S_i$ on $X_j$.
\item Fix $i\gg 0$ and an lc center $S_i$ of $(X_i,B_i,\Mm)$ such that $\phi_{i,j}$ induces an isomorphism for every lc center of $(S_i,B_{S_i},\Mm^{S_i})/U$ for any $j\geq i$. Then
\begin{enumerate}
    \item $\phi_{j,k}|_{S_{i,j}}: S_{i,j}\dashrightarrow S_{i,k}$ is an isomorphism in codimension $1$ for any $j,k\gg i$, and
    \item $B_{S_{i,j}}$ is the strict transform of $B_{S_{i,k}}$ for any $j,k\gg i$.
\end{enumerate}
\item Suppose that this $(K_X+B+\Mm_X)$-MMP$/U$ is a MMP with scaling of an $\Rr$-divisor $A\geq 0$ on $X$. Let $$\lambda_j:=\inf\{t\mid t\geq 0, K_{X_j}+B_j+tA_j+\Mm_{X_j}\text{ is nef/}U\}$$
be the scaling numbers, where $A_j$ is the strict transform of $A$ on $X_j$ for each $j$. Fix $i\gg 0$ and an lc center $S_i$ of $(X_i,B_i,\Mm)$ such that $\phi_{j,k}|_{S_{i,j}}: S_{i,j}\dashrightarrow S_{i,k}$ is an isomorphism in codimension $1$ and $B_{S_{i,j}}$ is the strict transform of $B_{S_{i,k}}$ for any $k,j\geq i$. Let $T$ be the normalization of the image of $S_i$ on $U$, $(S_{i}',B_{S_{i}'},\Mm^{S_i})$ a dlt model of $(S_{i},B_{S_{i}},\Mm^{S_i})$, and $A_{S_{i}'}$ the pullback of $A_i$ on $S_{i}'$. Then this $(K_X+B+\Mm_X)$-MMP$/U$ with scaling of $A$ induces the following commutative diagram$/T$
 \begin{center}$
 \xymatrixrowsep{0.135in}
\xymatrixcolsep{0.03in}
\xymatrix{
 (S_{i}',B_{S_{i}'},\Mm^{S_i})\ar@{-->}[rr]\ar@{->}[dd] & &  (S_{i,i+1}',B_{S_{i,i+1}'},\Mm^{S_i})\ar@{-->}[rr]\ar@{->}[dd] & &  \dots\ar@{-->}[rr] & & (S_{i,j}',B_{S_{i,j}'},\Mm^{S_i})\ar@{-->}[rr]\ar@{->}[dd] & & \dots\\
 && && && &&\\
 (S_{i},B_{S_{i}},\Mm^{S_i})\ar@{-->}[rr] & &  (S_{i,i+1},B_{S_{i,i+1}},\Mm^{S_i})\ar@{-->}[rr]& &  \dots\ar@{-->}[rr] & & (S_{i,j},B_{S_{i,j}},\Mm^{S_i})\ar@{-->}[rr] & & \dots\\
}$
\end{center}
such that
\begin{enumerate}
    \item $$(S_{i}',B_{S_{i}'},\Mm^{S_i})\dashrightarrow  (S_{i,i+1}',B_{S_{i,i+1}'},\Mm^{S_i})\dashrightarrow\dots\dashrightarrow (S_{i,j}',B_{S_{i,j}'},\Mm^{S_i})\dashrightarrow\dots$$
is a $(K_{S_{i}'}+B_{S_{i}'}+\Mm^{S_i}_{S_{i}'})$-MMP$/T$ with scaling of $A_{S_i'}$. Note that it is possible that $(S_{i,j}',B_{S_{i,j}'},\Mm^{S_i})\dashrightarrow (S_{i,j+1}',B_{S_{i,j+1}'},\Mm^{S_i})$ is the identity morphism or a composition of several steps of the  $(K_{S_{i,j}'}+B_{S_{i,j}'}+\Mm^{S_i}_{S_{i,j}'})$-MMP$/T$ for some $j$,
\item for any $j\geq i$, $(S_{i,j}',B_{S_{i,j}'},\Mm^{S_i})$ is a dlt model of $(S_{i,j},B_{S_{i,j}},\Mm^{S_i})$, and
\item let
$$\mu_j:=\inf\{t\mid t\geq 0, K_{S_{i,j}'}+B_{S_{i,j}'}+tA_{S_{i,j}'}+\Mm^{S_i}_{S_{i,j}'}\text{ is nef}/T\}$$
for each $j\geq i$, where $A_{S_{i,j}'}$ is the pullback of $A_j$ on $S_{i,j}'$. Then $\mu_j\leq\lambda_j$ for each $j\geq i$.
\end{enumerate}

\end{enumerate}
\end{lem}
\begin{proof}
Let $\Ii\subset [0,1]$ be a finite set such that $B\in\Ii$, and let $\Ii'\subset [0,+\infty)$ be a finite set such that $\Mm=\sum \mu_i\Mm_i$, where each $\Mm_i$ is nef$/U$ $\bb$-Cartier and each $\mu_i\in\Ii'$. Let $\phi_i:=\phi_{i,i+1}$ for each $i$.

We may assume that the MMP does not terminate, otherwise there is nothing left to prove. Possibly replacing $X$ with $X_i$ for $i\gg 0$, we may assume that each $\phi_i$ is a flip. Since the number of lc centers of $(X,B,\Mm)$ is finite, possibly replacing $X$ with $X_i$ for $i\gg 0$, we may assume that the flipping locus of $\phi_i$ does not contain any lc centers. This proves (1). 

We prove (2). We let $S:=S_i$. By (1), we may let $S_j:=S_{i,j}$ for any $j\geq i$. Possibly replacing $X$ we $X_i$, we may assume that $i=0$. By \cite[Proposition 2.10]{HL22}, for any $j$, the g-pair $(S_j,B_{S_j},\Mm^{S})$ given by the adjunction
$$K_{S_j}+B_{S_j}+\Mm_{S_j}^S:=(K_{X_j}+B_j+\Mm_{X_j})|_{S_j}$$
is dlt, and $B_{S_j}\in\mathbb S(\Ii,\Ii')$. By assumption, $\phi_{j,k}$ induces an isomorphism on $\lfloor B_{S_j}\rfloor$ for any $j,k$. Thus for any $j$ and any prime divisor $E$ over $S_j$, $\Center_{S_j}E\subset\lfloor B_{S_j}\rfloor$ if and only if $\Center_{S_{j+1}}E\subset\lfloor B_{S_{j+1}}\rfloor$. By the negativity lemma, $a(E,S_j,B_{S_j},\Mm^{S})\leq a(E,S_{j+1},B_{S_{j+1}},\Mm^{S})$ for each $j$ and any prime divisor $E$ over $S_j$. Thus
$$d_{\Ii,\Ii'}(S_j,B_{S_j},\Mm^{S})\geq d_{\Ii,\Ii'}(S_{j+1},B_{S_{j+1}},\Mm^{S})$$
for each $j$. Moreover, for any $j$ such that $S_j$ and $S_{j+1}$ are not isomorphic in codimension $1$,  if there exists a prime divisor $E$ on $S_{j+1}$ that is exceptional over $S_j$, then $$1-\alpha=a(E,S_{j+1},B_{S_{j+1}},\Mm^{S})>a(E,S_j,B_{S_j},\Mm^{S})$$ for some $\alpha\in\mathbb S(\Ii,\Ii')$, and hence $$d_{\Ii,\Ii'}(S_j,B_{S_j},\Mm^{S})> d_{\Ii,\Ii'}(S_{j+1},B_{S_{j+1}},\Mm^{S}).$$
By \cite[Remark 4.6]{HL22}, $d_{\Ii,\Ii'}(S_j,B_{S_j},\Mm^{S})<+\infty$. Thus possibly replacing $X$ with $X_j$ for some $j\gg 0$, we may assume that
$d_{\Ii,\Ii'}(S_j,B_{S_j},\Mm^{S})=d_{\Ii,\Ii'}(S_k,B_{S_k},\Mm^{S})$ for any $j,k$. Thus $S_j\dashrightarrow S_{j+1}$ does not extract any divisor for any $j$. In particular, $\rho(S_{j+1})\leq \rho(S_j)$, and $\rho(S_{j+1})<\rho(S_j)$ if $S_j\dashrightarrow S_{j+1}$ contracts a divisor. Thus possibly replacing $X$ with $X_j$ for some $i\gg 0$, we may assume that  $S_j$ and $S_{j+1}$ are isomorphic in codimension $1$ for each $j$, which implies (2.a). Since $d_{\Ii,\Ii'}(S_j,B_{S_j},\Mm^{S})=d_{\Ii,\Ii'}(S_k,B_{S_k},\Mm^{S})$ for any $j,k$, (2.b) follows from (2.a).

We prove (3). Since $i\gg 0$, possibly replacing $X$ with $X_i$, we may assume that $i=0$ and $\phi_{j}$ is a flip for every $j$. We let $S:=S_0,S':=S_0',S_j:=S_{0,j}$, and $S_j':=S_{0,j}'$ for every $j$. We let $X_j\rightarrow Z_j\leftarrow X_{j+1}$ be each flip and let $T_j$ be the normalization of the image of $S_{j}$ on $Z_j$ for each $j$. Then we have an induced birational map $S_{j}\dashrightarrow S_{j+1}$ for each $j$.

Since $\phi_{0}$ is a $(K_{X_0}+B_0+\Mm_{X_0})$-flip$/U$, $X_{1}\rightarrow Z_0$ is $(K_{X_{1}}+B_1+\Mm_{X_1})$-positive and $K_{S_1}+B_{S_1}+\Mm^S_{S_1}$ is ample$/T_0$. In particular, $(S_1,B_{S_1},\Mm^S)/T_0$ is a weak lc model of $(S_0,B_{S_0},\Mm^S)$. By \cite[Lemmas 3.9, 3.15]{HL21a} and \cite[Theorem 3.14]{HL21a}, we may run a $(K_{S_0'}+B_{S_0'}+\Mm^{S}_{S_0'})$-MMP$/T_0$ with scaling of an ample$/T_0$ divisor, which terminates with a good minimal model of $(S_0',B_{S_0'},\Mm^S)/T_0$. By \cite[Lemma 3.9]{HL21a}, $(S_0',B_{S_0'},\Mm^S)$ is a dlt model of $(S_1,B_{S_1},\Mm^S)$. Since 
$$K_{S_0'}+B_{S_0'}+\lambda_0A_{S_0'}+\Mm^{S}_{S_0'}\equiv_{T_0}0,$$
 this MMP is also a $(K_{S_0'}+B_{S_0'}+\Mm^{S}_{S_0'})$-MMP$/T_0$ with scaling of $\lambda_0A_{S_0'}$. 
 We may replace $(S_0,B_{S_0},\Mm^S)/T$ with $(S_1,B_{S_1},\Mm^S)/T$ and continue this process. This gives us the desired $(K_{S_0'}+B_{S_0'}+\Mm^{S}_{S_0'})$-MMP$/T$ with scaling of $A_{S_0'}$, which gives the commutative diagram, and proves (3.a) and (3.b). For each $j$, since $K_{S_j'}+B_{S_j'}+\lambda_jA_{S_j'}+\Mm^{S}_{S_j'}\equiv_{T_j}0$, 
 $K_{S_j'}+B_{S_j'}+\lambda_jA_{S_j'}+\Mm^{S}_{S_j'}$ is nef, hence $\mu_j\leq\lambda_j$, and we get (3.c).
\end{proof}

\section{Relative Nakayama-Zariski decomposition}\label{sec: relative nz decomposition}

\begin{defn}\label{defn: rel nz decomposition}
Let $\pi: X\rightarrow U$ be a projective morphism from a normal variety to a variety, $A$ an ample$/U$ $\Rr$-divisor on $X$, $D$ a pseudo-effective$/U$ $\Rr$-Cartier $\Rr$-divisor on $X$, and $P$ a prime divisor on $X$. For any big$/U$ $\Rr$-Cartier $\Rr$-divisor $B$, we define
$$\sigma_P(X/U,B):=\inf\{\mult_PB'\mid 0\leq B'\sim_{\Rr,U}B\}.$$
We define
$$\sigma_P(X/U,D):=\lim_{\epsilon\rightarrow0^+}\sigma_P(X/U,D+\epsilon A),$$
where we allow $+\infty$ as a limit as well. As in \cite[III \S 1.]{Nak04}, we can easily check that $\sigma_P(X/U,D)$ is well-defined and does not depend on the choice of $A$ (we left the proof for the readers). We let
$$N_{\sigma}(X/U,D):=\sum_{C\text{ is a prime divisor on }X}\sigma_C(X/U,D)\cdot C$$ 
be a formal sum of divisors with coefficients in $\Rr_{\ge0}\cup\{+\infty\}$. We say that $D$ is movable$/U$ if $N_{\sigma}(X/U,D)=0$, and this coincides with the original definition when $D$ is big$/U$. 

For any divisor $D'$ on $X$, we say $D'\le N_{\sigma}(X/U,D)$ if $\mult_C(D')\le \sigma_C(X/U,D)$ for any prime divisor $C$ on $X$. We can naturally define the addition of $D'$ and $N_{\sigma}(X/U,D)$ as
$$
N_{\sigma}(X/U,D)+D':=\sum_{C\text{ is a prime divisor on }X}(\sigma_C(X/U,D)+\mult_{C}(D'))\cdot C,
$$
by noticing that $+\infty+a=+\infty$ for any $a\in\Rr$. If $f:X\to Y$ is a projective birational morphism over $U$, then we can define the pushforward
$$
f_*N_{\sigma}(X/U,D):=\sum_{C\text{ is a prime divisor on }X}\sigma_C(X/U,D)\cdot f_*C
$$
as a formal sum of divisors with coefficients in $\mathbb R_{\geq 0}\cup\{+\infty\}$.

We define the support of $N_{\sigma}(X/U,D)$ as 
$$
\Supp N_{\sigma}(X/U,D):=\bigcup_{\sigma_C(X/U,D)>0}C.
$$

If there are only finitely many prime divisors $C$ on $X$ such that $\sigma_C(X/U,D)>0$ and $\sigma_C(X/U,D)<+\infty$ (e.g $D\ge0$), then (we say) $N_{\sigma}(X/U,D)$ is well-defined as a divisor and we let
$$P_{\sigma}(X/U,D):=D-N_{\sigma}(X/U,D).$$
\end{defn}

Definition \ref{defn: rel nz decomposition} is the same as the one adopted in \cite{HX13,HMX18}. The following lemma shows that the relative Nakayama-Zariski decomposition defined in Definition \ref{defn: rel nz decomposition} is the same as the $\sigma$-decomposition defined in \cite[III. \S 4.a]{Nak04}:

\begin{lem}\label{lem: nz decomposition definitions are equivalent}
Notation as in Definition \ref{defn: rel nz decomposition}. If $X$ is smooth, then $\sigma_P(X/U,D)$ is the same as $\sigma_P(D,X/U)$, where the latter is the value defined as in Nakayama's original relative $\sigma$-decomposition \cite[III. \S 4.a]{Nak04}. 
\end{lem}
\begin{proof}
By definition, we only need to deal with the case when $D$ is big. We may pick an affine open subset $U^0$ of $U$ such that $P$ intersects $X^0:=X\times_UU^0$. Let $P^0:=P\times_UU^0$ and $D^0:=D\times_UU^0$. Then
$$\sigma_P(X/U,D)=\sigma_{P^0}(X^0/U^0,D^0).$$
Possibly replacing $(X/U,D)$ and $P$ with $(X^0/U^0,D^0)$ and $P^0$ respectively, we may assume that $U$ is affine. Thus for any Cartier divisor $Q$ on $U$, there exists a principal divisor $Q'$ on $U$ such that $Q'=Q$ in a neighborhood of the generic point of $\pi(P)$. In particular, we have
$$\sigma_P(X/U,D)=\inf\{\mult_{P^0}B'\mid 0\leq B'\sim_{\Rr}B^0\}.$$
For any Cartier divisor $F$ on $X$, let
$$m_F:=\inf\{+\infty,\mult_PF'\mid 0\leq F'\sim F\}.$$
If $m_F<+\infty$, then by definition,
$$m_F=\max\{m\in\mathbb N\mid H^0(X,F-mP)\hookrightarrow H^0(X,F)\text{ is an isomorphism}\}.$$
Moreover, since $U$ is affine and $H^0(X,\mathcal{O}_X(F))=H^0(U,\pi_*\mathcal{O}_X(F))$, if $m_F<+\infty$, then
$$m_F=\max\{m\in\mathbb N\mid \pi_*\mathcal{O}_X(F-mP)\hookrightarrow\pi_*\mathcal{O}_X(F)\text{ is an isomorphism}\}.$$
Now the lemma follows from the construction in \cite[III. \S 4.a]{Nak04}.
\end{proof}

\begin{lem}\label{lem: nz pre properties}
Let $\pi: X\rightarrow U$ be a projective morphism from a normal variety to a variety, $D,D'$ two pseudo-effective$/U$ $\Rr$-Cartier $\Rr$-divisors on $X$, and $P$ a prime divisor on $X$. 
\begin{enumerate}
    \item If $D$ is nef$/U$, then $\sigma_P(X/U,D)=0$. 
    \item $\sigma_P(X/U,D+D')\leq\sigma_P(X/U,D)+\sigma_P(X/U,D')$.
    \item If $\sigma_P(X/U,D')<+\infty$, then $\lim_{\epsilon\rightarrow 0^+}\sigma_P(X/U,D+\epsilon D')=\sigma_P(X/U,D)$.
\end{enumerate}
\end{lem}
\begin{proof}
Let $A$ be an ample$/U$ divisor on $X$. \par
(1) is straightforward from the definition.

(2) follows from the fact that $\sigma_P(X/U,D+D'+\epsilon A)\le\sigma_P(X/U,D+\frac{\epsilon}{2}A)+\sigma_P(X/U, D'+\frac{\epsilon}{2}A)$.

There exists $a>0$ such that $A-aD'$ is ample$/U$. Thus, by (1) and (2), we have
$$\sigma_P(X/U,D)+\sigma_P(X/U,a\epsilon D')\ge\sigma_P(X/U, D+a\epsilon D')\ge\sigma_P(X/U,D+\epsilon A),$$
and (3) follows after taking $\epsilon\to 0^+$.
\end{proof}

\begin{lem}\label{lem: nz keep under pullback}
Let $\pi: X\rightarrow U$ be a projective morphism from a normal variety to a variety and $D$ a pseudo-effective$/U$ $\Rr$-Cartier $\Rr$-divisor on $X$. Let $f: Y\rightarrow X$ be a projective birational morphism. Then:
\begin{enumerate}
    \item For any prime divisor $P$ on $X$, we have
$$\sigma_P(X/U,D)=\sigma_{f^{-1}_*P}(Y/U,f^*D).$$
    \item For any exceptional$/X$ $\Rr$-Cartier $\Rr$-divisor $E\ge0$ and any prime divisor $P$ on $Y$, we have $$\sigma_P(Y/U,f^*D+E)=\sigma_P(Y/U,f^*D)+\mult_PE.$$
    \item For any exceptional$/X$ $\Rr$-Cartier $\Rr$-divisor $E\ge0$ on $Y$ we have 
    $$N_{\sigma}(X/U,D)=f_*N_{\sigma}(Y/U,f^*D+E)$$
    as a formal sum of divisors with coefficients in $\mathbb R_{\geq 0}\cup\{+\infty\}$. In particular, if $N_{\sigma}(Y/U,f^*D+E)$ is well-defined, then
    $N_{\sigma}(X/U,D)$ is well-defined.
    \item If $D'\ge0$ is an $\Rr$-Cartier $\Rr$-divisor on $X$ such that $D'\le N_\sigma(X/U,D)$, then $f^*D'\le N_\sigma(Y/U,f^*D)$ and $D-D'$ is pseudo-effective$/U$. 
\end{enumerate}
\end{lem}
\begin{proof}
Set $g=\pi\circ f$ and let $A$ (resp. $A'$) be an ample$/U$ divisor on $X$ (resp. $Y$). Fix a real number $a>0$ such that $aA'+f^*A$ is ample$/U$. Notice that $f^{-1}_*P$ is a prime divisor on $Y$. Since $f^*A$ is semi-ample$/U$, by Lemma \ref{lem: nz pre properties} we have $\lim_{\epsilon\to0^+}\sigma_{f^{-1}_*P}(Y/U,f^*D+\epsilon f^*A)=\sigma_{f^{-1}_*P}(Y/U,f^*D)$.

Since $\pi_*\Oo_X(F)=g_*\Oo_Y(f^*F)$ for any Cartier divisor $F$ on $X$, by definition we have $\sigma_P(X/U,D+\epsilon A)=\sigma_{f_*^{-1}P}(Y/U,f^*D+\epsilon f^*A)$ for any $\epsilon>0$. Thus we have 
$$
\sigma_P(X/U,D)=\lim_{\epsilon\to0^+}\sigma_{f_*^{-1}P}(Y/U,f^*D+\epsilon f^*A)=\sigma_{f^{-1}_*P}(Y/U,f^*D)
$$ 
which is (1).

Since $\lim_{\epsilon\to0^+}\sigma_P(Y/U,f^*D+\epsilon f^*A)=\sigma_P(Y/U,f^*D)$, we may assume that $D$ is a big$/U$. (2) follows from the fact that $g_*\Oo_Y(f^*F+E)=\pi_*\Oo_X(F)$ for any Cartier divisor $F$ on $X$ and any exceptional$/X$ divisor $E\ge0$.

We have 
$$N_{\sigma}(Y/U,f^*D+E)=N_{\sigma}(Y/U,f^*D)+E$$
by (2) and
$$f_*N_{\sigma}(Y/U,f^*D)=N_{\sigma}(X/U,D)$$
by (1), which imply (3).

For (4), since there are only finitely many prime divisors $P$ on $X$ such that $\mult_PD'>0$, by assumption and by the definition of $\sigma_P(X/U,D)$ we know that $D'\le D_{\epsilon}''$ for any element $D''_\epsilon\sim_{\Rr,U}D+\epsilon A$ and any $1\gg\epsilon>0$. Then $[D-D']=\lim_{\epsilon\to0^+}[D''_\epsilon-D']$ is indeed pseudo-effective$/U$. Moreover, $f^*D'\le f^*D''_\epsilon$ for any $1\gg\epsilon>0$ and Lemma \ref{lem: nz pre properties}(3) implies that $\mult_Pf^*D'\le\sigma_P(Y/U,f^*D)$ for any prime divisor $P$ on $Y$ by the same argument as in the proof of (2) above.
\end{proof}

\begin{lem}\label{lem: nz finite is well-defined}
Let $\pi: X\rightarrow U$ be a projective morphism from a normal variety to a variety and $D$ a pseudo-effective$/U$ $\Rr$-Cartier $\Rr$-divisor on $X$. Then there are only finitely many prime divisors $P$ on $X$ such that $\sigma_P(X/U,D)>0$. In particular, $\Supp N_\sigma(X/U,D)$ can be regarded as a reduced divisor. If furthermore $\sigma_P(X/U,D)<+\infty$ for any prime divisor $P$ on $X$, then $N_{\sigma}(X/U,D)$ and $P_{\sigma}(X/U,D)$ are well-defined as divisors.
\end{lem}
\begin{proof}
Let $f: Y\rightarrow X$ be a resolution of $X$. By Lemma \ref{lem: nz keep under pullback}(1), for any prime divisor $P$ on $X$ such that $\sigma_{P}(X/U,D)\not=0$, $\sigma_{f^{-1}_*P}(Y/U,f^*D)\not=0$. Therefore, we only need to show that there are finitely many prime divisors $P_Y$ on $Y$ such that $\sigma_{P_Y}(Y/U,f^*D)\not=0$. Possibly replacing $X$ with $Y$ and $D$ with $f^*D$, we may assume that $X$ is smooth. In the following, we will show that there are at most $\rho(X/U)=\dim N^1(X/U)_{\Rr}$ prime divisors $P$ on $X$ such that $\sigma_P(X/U,D)\not=0$. 

Let $P_1,P_2,...,P_l$ be distinct prime divisors of $X$ such that $\sigma_{P_i}(X/U,D)>0$ for each $i$.  If $l\leq\dim N^1(X/U)_{\Rr}$ then we are done. Otherwise, $P_1,P_2,...,P_l$ are not linearly independent in $N^1(X/U)_\Rr$ and possibly reordering indices, we have
$$\sum_{i=1}^{s}x_iP_i\equiv_U\sum_{j=s+1}^{l}x_jP_j\in N^1(X/U)$$
for some $x_1,x_2,...,x_l\in\Rr_{\ge0}$ and $1\le s\le l$, we may also assume that $x_1\neq 0$. By Lemma \ref{lem: nz decomposition definitions are equivalent} and \cite[III, Lemma 4.2(2)]{Nak04}, we have
$$\sigma_{P_i}(X/U,\sum^{l}_{j=1}x_jP_j)=x_i$$
for any $x_1,x_2,...,x_l\in\Rr_{\ge0}$. Since $\sigma_P(X/U,D)$ depends only on the numerical equivalence class of $D$ over $U$, by Lemma \ref{lem: nz decomposition definitions are equivalent} and \cite[III, Lemma 4.2(2)]{Nak04} again, we obtain
$$
x_1=\sigma_{P_1}(X/U,\sum^{s}_{i=1}x_iP_i)=\sigma_{P_1}(X/U,\sum^{l}_{j=s+1}x_jP_j)=0,
$$
a contradiction.
\end{proof}

\begin{defn}\label{defn: sigma over X}
Let $\pi: X\rightarrow U$ be a projective morphism from a normal variety to a variety, $D$ a pseudo-effective$/U$ $\Rr$-Cartier $\Rr$-divisor on $X$, and $P$ a prime divisor over $X$. Let $f: Y\rightarrow X$ be a projective birational morphism such that $P$ descends to $Y$. We define
$$\sigma_P(X/U,D):=\sigma_P(Y/U,f^*D).$$
By Lemma \ref{lem: nz keep under pullback}, $\sigma_P(X/U,D)$ is independent of the choice of $Y$. Also notice that $f^*D$ is pseudo-effective$/U$ iff $D$ is.
\end{defn}

\begin{lem}\label{lem: nz basic properties}
Let $\pi: X\rightarrow U$ be a projective morphism from a normal variety to a variety, $D,D'$ two pseudo-effective$/U$ $\Rr$-Cartier $\Rr$-divisors on $X$, and $P$ a prime divisor over $X$. Then
\begin{enumerate}
    \item $\sigma_P(X/U,D+D')\leq\sigma_P(X/U,D)+\sigma_P(X/U,D')$.
    \item If $\sigma_P(X/U,D')<+\infty$, then $\lim_{\epsilon\rightarrow 0^+}\sigma_P(X/U,D+\epsilon D')=\sigma_P(X/U,D)$.
    \item If $D$ is a movable$/U$ $\Rr$-Cartier $\Rr$-divisor, then $N_{\sigma}(X/U,D)=0$ and $P_{\sigma}(X/U,D)=D$ is movable.
    \item $\Supp N_{\sigma}(X/U,D)$ coincides with the divisorial part of ${\bf{B}}_{-}(D/U)$.
    \item If $0\leq D'\leq N_{\sigma}(X/U,D)$, then $N_{\sigma}(X/U,D-D')+D'=N_{\sigma}(X/U,D)$.
    \item If $D'\geq 0$ and $\Supp D'\subset\Supp N_{\sigma}(X/U,D)$, then $N_{\sigma}(X/U,D+D')=N_{\sigma}(X/U,D)+D'$.
\end{enumerate}
\end{lem}
\begin{proof}
Let $A$ be an ample$/U$ $\Rr$-divisor on $X$. 

(1) and (2) follow directly from Lemma \ref{lem: nz pre properties}(2)(3).



For (3), if this is not true, then we have $\sigma_P(X/U,D)>0$ for some $P$. By definition, there exist an $\epsilon>0$ such that $\sigma_P(X/U,D+\epsilon A)>0$. Assume $[D]=\lim_{i\to\infty}[D_i]$, where $D_i$ is a movable divisor for each $i\ge1$. Then $\epsilon A-(D_i-D)$ is ample for any $i\gg0$, and we have
$$
0<\sigma_P(X/U,D+\epsilon A)=\sigma_P(X/U,D_i+\epsilon A-(D_i-D))\le\sigma_P(X/U,D_i)=0,
$$ which is a contradiction.

For (4), from the definition of ${\bf{B}}_-(D/U)$ we know that $\Supp N_\sigma(X/U,D)\subset{\bf{B}}_-(D/U)$. For any divisorial component $P$ of ${\bf{B}}_-(D/U)$, there exist $\epsilon>0$ such that $P\subset {\bf{B}}(D+\epsilon A/U)$, so $\sigma_P(X/U,D)\ge\sigma_P(X/U,D+\epsilon A)>0$.

By Lemma \ref{lem: nz keep under pullback}, possibly replacing $X$ with a resolution, we may assume that $X$ is smooth and $P$ is a prime divisor on $X$. Then (5) and (6) follow from Lemma \ref{lem: nz decomposition definitions are equivalent} and \cite[III, Lemma 4.2]{Nak04}. Notice that $N_\sigma(X/U,D-D')$ makes sense by Lemma \ref{lem: nz keep under pullback}(4).
\end{proof}

\begin{lem}[{cf. \cite[Lemma 2.4]{Has20}}]\label{lem: Has20 2.4 rel ver}
Let $\pi: X\rightarrow U$ be a projective morphism from a normal variety to a variety, $D$ (resp. $D'$) a pseudo-effective$/U$ $\Rr$-Cartier $\Rr$-divisor on $X$ such that $N_{\sigma}(X/U,D)$ (resp. $N_\sigma(X/U,D')$) is well-defined as a divisor. Then there exists $t_0>0$ such that $\Supp N_{\sigma}(X/U,D+tD')$ is independent of $t$ for any $t\in (0,t_0]$.
\end{lem}
\begin{proof}
The proof is exactly the same as that of \cite[Lemma 2.4]{Has20}.
\end{proof}

\begin{lem}\label{lem: nz for lc divisor}
Let $(X,B,\Mm)/U$ be a $\Qq$-factorial NQC dlt g-pair. Then for any partial $(K_X+B+\Mm_X)$-MMP$/U$ $\phi: X\dashrightarrow \bar X$, 
    \begin{enumerate}
    \item the divisors contracted by $\phi$ are contained in $\Supp N_{\sigma}(X/U,K_X+B+\Mm_X)$, and
    \item let $\bar B$ be the strict transform of $B$ on $\bar X$. If $K_{\bar X}+\bar B+\Mm_{\bar X}$ is movable$/U$, then $\Supp N_{\sigma}(X/U,K_X+B+\Mm_X)$ is the set of all $\phi$-exceptional divisors.
    \end{enumerate}
 \end{lem}
\begin{proof}
Let $p: W\rightarrow X$ and $q: W\rightarrow \bar X$ be a resolution of indeterminacies of $\phi$. Then
$$p^*(K_X+B+\Mm_X)=q^*(K_{\bar X}+\bar B+\Mm_{\bar X})+E$$
for some $E\geq 0$ that is exceptional$/\bar X$ and $\Supp E$ contains the strict transforms on $W$ of all $\phi$-exceptional divisors. By Lemma \ref{lem: nz keep under pullback}(2) we have
$$
\Supp E\subset\Supp N_\sigma(W/U,q^*(K_{\bar X}+\bar B+\Mm_{\bar X})+E)=\Supp N_\sigma(W/U,p^*(K_X+B+\Mm_X)),
$$
and by Lemma \ref{lem: nz keep under pullback}(3) we know that $\Supp p_*E\subset\Supp N_\sigma(X/U,K_X+B+\Mm_X)$. Therefore, any $\phi$-exceptional divisor is contained in $\Supp N_\sigma(X/U,K_X+B+\Mm_X)$.

If $K_{\bar X}+\bar B+\Mm_{\bar X}$ is movable$/U$, then by Lemma \ref{lem: nz keep under pullback}(3) we have $q_* N_\sigma(W/U,q^*(K_{\bar X}+\bar B+\Mm_{\bar X})+E)=0$ so $\Supp N_\sigma(W/U,q^*(K_{\bar X}+\bar B+\Mm_{\bar X})+E)$ (viewed as a reduced divisor) is $q$-exceptional. By Lemma \ref{lem: nz keep under pullback}(3) again we have $\Supp N_\sigma(X/U,K_X+B+\Mm_X)=\Supp p_*N_\sigma(W/U,q^*(K_{\bar X}+\bar B+\Mm_{\bar X})+E)$, whose components are all $\phi$-exceptional.
\end{proof}

\begin{lem}\label{lem: hmx18 2.7.3 gpair rel}
Let $(X,B,\Mm)/U$ be an NQC lc g-pair such that $K_X+B+\Mm_X$ is pseudo-effective$/U$. Let $\phi: X\dashrightarrow X'$ be a birational map$/U$ which does not extract any divisor and $B'$ the strict transform of $B$ on $X'$, such that
\begin{enumerate}
\item $K_{X'}+B'+\Mm_{X'}$ is nef$/U$, and
\item $\phi$ only contracts divisors contained in $\Supp N_{\sigma}(X/U,K_X+B+\Mm_X)$.
\end{enumerate}
Then $(X',B',\Mm)/U$ is a log minimal model (not necessarily $\Qq$-factorial) of $(X,B,\Mm)/U$.
\end{lem}
\begin{proof}
Let $p: W\rightarrow X$ and $q: W\rightarrow X'$ be a resolution of indeterminacies of $\phi$ such that
$$p^*(K_X+B+\Mm_X)+E=q^*(K_{X'}+B'+\Mm_{X'})+F$$
where $E\geq 0,F\geq 0$, and $E\wedge F=0$. Then $E$ and $F$ are $q$-exceptional. By Lemma \ref{lem: nz basic properties}(3) and \ref{lem: nz keep under pullback}(2)(3), $F=N_{\sigma}(W/U,q^*(K_{X'}+B'+\Mm_{X'})+F)$. 

We may write $E=E_1+E_2$ such that $E_1$ is $p$-exceptional and every component of $E_2$ is not $p$-exceptional. Then $p_*E_2$ is $\phi$-exceptional and therefore by assumption (2), we obtain $\Supp p_*E_2\subset\Supp N_{\sigma}(X/U,K_X+B+\Mm_X)$. By Lemma \ref{lem: nz keep under pullback}(3), we know 
$$
\Supp E_2\subset\Supp N_{\sigma}(W/U,p^*(K_X+B+\Mm_X)).
$$
By Lemma \ref{lem: nz keep under pullback}(2), we have
$$
N_\sigma(W/U,p^*(K_X+B+\Mm_X)+E_1)=N_\sigma(W/U,p^*(K_X+B+\Mm_X))+E_1.
$$
Therefore, $$\Supp(E_1+E_2)\subset\Supp N_{\sigma}(W/U,p^*(K_X+B+\Mm_X)+E_1),$$ 
and then by Lemma \ref{lem: nz basic properties}(6), we have 
\begin{align*}
    N_\sigma(W/U,p^*(K_X+B+\Mm_X)+E_1+E_2)&=N_\sigma(W/U,p^*(K_X+B+\Mm_X)+E_1)+E_2\\
    &=N_\sigma(W/U,p^*(K_X+B+\Mm_X))+E_1+E_2,
\end{align*}
which immediately implies that
$$
\Supp E=\Supp(E_1+E_2)\subset\Supp N_{\sigma}(W/U,p^*(K_X+B+\Mm_X)+E_1+E_2)=\Supp F,
$$
and hence $E$ must be zero. Now by Lemma \ref{lem: nz keep under pullback}(3) again $\Supp p_*F=\Supp N_\sigma(X/U,K_X+B+\Mm_X)$, which contains all $\phi$-exceptional divisors and we are done.
\end{proof}

\section{Reduction via Iitaka fibration}

\begin{lem}\label{lem: has19 3.2 step 3 abu ver}
Let $(X,B,\Mm)/U$ be a $\Qq$-factorial NQC lc g-pair with $X$ klt and $\pi: X\rightarrow U$ the induced morphism, such that
\begin{enumerate}
\item $\pi$ is an equidimensional contraction,
\item $U$ is quasi-projective and $\Qq$-factorial, and
\item $\kappa_{\sigma}(X/U,K_X+B+\Mm_X)=\kappa_{\iota}(X/U,K_X+B+\Mm_X)=0$.
\end{enumerate}
Let $A\geq 0$ be an ample$/U$ $\Rr$-divisor on $X$ such that $(X,B+A,\Mm)$ is lc and $K_X+B+A+\Mm_X$ is nef$/U$, and run a $(K_X+B+\Mm_X)$-MMP$/U$ with scaling of $A$. Then this MMP terminates with a good minimal model $(X',B',\Mm)/U$ of $(X,B,\Mm)/U$. Moreover, $K_{X'}+B'+\Mm_{X'}\sim_{\Rr,U} 0$.
\end{lem}
\begin{proof}
If $\dim X=\dim U$, then $\pi$ is the identity map since $\pi$ is an equidimensional contraction and there is nothing left to prove. In the following, we assume that $\dim X>\dim U$.

Since $\kappa_{\iota}(X/U,K_X+B+\Mm_X)=0$, $K_X+B+\Mm_X\sim_{\Rr,U}E\geq 0$ for some $\Rr$-divisor $E$ on $X$. We may write $E=E^{h}+E^{v}$, such that $E^h\geq 0,E^v\geq 0$, each component of $E^h$ is horizontal over $U$, and $E^v$ is vertical over $U$. Since $\pi$ is equidimensional, the image of any component of $E^v$ on $U$ is a divisor. Since $U$ is $\Qq$-factorial, for any prime divisor $P$ on $U$, we may define
$$\nu_P:=\sup\{\nu\geq 0 \mid E^v-\nu\pi^*P\geq 0\}.$$
Then $\nu_P>0$  for only finitely many prime divisors $P$ on $U$. Possibly replacing $E^v$ with $E^v-\pi^*(\sum_P\nu_PP)$, we may assume that $E^v$ is very exceptional over $U$. 

Let $F$ be a very general fiber of $\pi$. 
Let 
$$(X,B,\Mm):=(X_0,B_0,\Mm)\dashrightarrow (X_1,B_1,\Mm)\dashrightarrow\dots\dashrightarrow (X_i,B_i,\Mm)\dashrightarrow\dots$$
be a $(K_X+B+\Mm_X)$-MMP$/U$ with scaling of $A$, and let $A_i,E^h_i,E^v_i,F_i$ be the strict transforms of $A,E^h,E^v,F$ on $X_i$ respectively. Then we have
$$
\kappa_{\sigma}(X_i/U,K_{X_i}+B_i+\Mm_{X_i})=\kappa_{\iota}(X_i/U,K_{X_i}+B_i+\Mm_{X_i})=0
$$ by Lemma \ref{lem: property of numerical and Iitaka dimension}(6) and hence
$$
\kappa_\sigma(E^h_i|_{F_i})=\kappa_{\sigma}((K_{X_i}+B_{X_i}+\Mm_{X_i})|_{F_i})=\kappa_{\iota}((K_{X_i}+B_{X_i}+\Mm_{X_i})|_{F_i})=\kappa_{\iota}(E^h_i|_{F_i})=0,
$$
since $E_i^h|_{F_i}=(E_i^h+E_i^v)|_{F_i}\sim_{\Rr} (K_{X_i}+B_i+\Mm_{X_i})|_{F_i}$. As in the proof of \cite[Theorem 3.4]{Bir12}, there exists a positive integer $n$ such that $K_{X_n}+B_n+\Mm_{X_n}$ is a movable$/U$ $\Rr$-Cartier $\Rr$-divisor. Therefore the restriction $(K_{X_n}+B_n+\Mm_{X_n})|_{F_n}\sim_{\Rr}E^h_n|_{F_n}$ is also a movable $\Rr$-Cartier $\Rr$-divisor since $F$ is a very general fiber. In particular, $N_\sigma(E^h_n|_{F_n})=0$ by Lemma \ref{lem: nz basic properties}(3).
Notice that now $F_n$ is a normal projective variety and let $g: F'_n\to F_n$ be a resolution of singularities. Then $\kappa_\sigma(g^*(E^h_n|_{F_n}))=\kappa_\iota(g^*(E^h_n|_{F_n}))=0$. By \cite[V, 1.12 Corollary]{Nak04} we have 
$$
N_\sigma(g^*(E^h_n|_{F_n}))\equiv g^*(E^h_n|_{F_n}).
$$
Since $g^*(E^h_n|_{F_n})\ge0$, we have $N_\sigma(g^*(E^h_n|_{F_n}))\le g^*(E^h_n|_{F_n})$ by the definition. Hence we must have
$$
N_\sigma(g^*(E^h_n|_{F_n}))= g^*(E^h_n|_{F_n}).
$$
Therefore by Lemma \ref{lem: nz keep under pullback}(3) we get 
$$
E^h_n|_{F_n}=g_*(g^*(E^h_n|_{F_n}))=g_*N_\sigma(g^*(E^h_n|_{F_n}))=N_\sigma(E^h_n|_{F_n})=0.
$$
This immediately implies that $E^h_n=0$ since $E^h_n\geq 0$ is horizontal over $U$. Notice that our $(K_X+B+\Mm_X)$-MMP$/U$ is also a $(E^h+E^v)$-MMP$/U$ and $E^v$ is very exceptional over $U$. By Lemma \ref{lem: rlinear version of HL22 3.8}, this MMP terminates with a log minimal model $(X',B',\Mm)/U=(X_m,B_m,\Mm)/U$ of $(X,B,\Mm)/U$ for some positive integer $m$, such that $K_{X'}+B'+\Mm_{X'}\sim_{\Rr,U}0$. In particular, $(X',B',\Mm)/U$ is a good minimal model of $(X,B,\Mm)/U$.
\end{proof}

\begin{thm}\label{thm: HL21a v2 4.1 abu ver}
Let $(X,B,\Mm)/U$ be an NQC lc g-pair and $\pi: X\rightarrow V$ a contraction over $U$ such that $V$ is quasi-projective. Assume that $\kappa_{\sigma}(X/V,K_X+B+\Mm_X)=\kappa_{\iota}(X/V,K_X+B+\Mm_X)=0$. Then there exists a $\Qq$-factorial NQC dlt g-pair $(X',B',\Mm)/U$, a contraction $\pi': X'\rightarrow V'$ over $U$, and a birational projective morphism $\varphi: V'\rightarrow V$ over $U$ satisfying the following:
\begin{center}$\xymatrix{
X'\ar@{->}[d]_{\pi'}\ar@{-->}[rr]& & X\ar@{->}[d]^{\pi}\\
V'\ar@{->}[rr]^{\varphi}\ar@{->}[dr] & & V\ar@{->}[dl]\\
& U &
}$
\end{center}
\begin{enumerate}
    \item $X'$ is birational to $X$ and $V'$ is smooth,
    \item $K_{X'}+B'+\Mm_{X'}\sim_{\Rr,V'}0$.
    \item $(X,B,\Mm)/U$ has a good minimal model if and only if $(X',B',\Mm)/U$ has a good minimal model.
    \item Any weak lc model of $(X,B,\Mm)/U$ is a weak lc model of $(X',B',\Mm)/U$, and any weak lc model of $(X',B',\Mm)/U$ is a weak lc model of $(X,B,\Mm)/U$.
    \item If all lc centers of $(X,B,\Mm)$ dominate $V$, then all lc centers of $(X',B',\Mm)$ dominate $V'$.
    \item $\kappa_{\sigma}(X/U,K_X+B+\Mm_X)=\kappa_{\sigma}(X'/U,K_{X'}+B'+\Mm_{X'})$ and $\kappa_{\iota}(X/U,K_X+B+\Mm_X)=\kappa_{\iota}(X'/U,K_{X'}+B'+\Mm_{X'})$
\end{enumerate}
\end{thm}

\begin{proof} Let $h: W\rightarrow X$ be a log resolution of $(X,\Supp B)$ such that $\Mm$ descends to $W$. 
By \cite[Lemma 3.6]{HL21a}, $(X,B,\Mm)$ has a proper log smooth model $(W,B_W,\Mm)$ for some $\Rr$-divisor $B_W$ on $W$. By Lemmas \ref{lem: property of numerical and Iitaka dimension}(3) and \cite[Lemma 3.7]{HL21a}, \cite[Theorem 3.14]{HL21a}, and \cite[Lemmas 3.10, 3.17]{HL21a}, we may replace $(X,B,\Mm)$ with $(W,B_W,\Mm)$, and assume that $(X,B)$ is log smooth and $\Mm$ descends to $X$. 

By Theorem \ref{thm: has19 weak semistable reduction}, there exists a commutative diagram of projective morphisms
\begin{center}$\xymatrix{
Y\ar@{->}[r]^{f}\ar@{->}[d]_{\pi_Y} & X\ar@{->}[d]^{\pi}\\
V'\ar@{->}[r]^{\varphi} & V
}$
\end{center}
such that
\begin{itemize}
    \item $f,\varphi$ are birational morphisms, $\pi_Y$ is an equidimensional contraction, $Y$ only has $\Qq$-factorial toroidal singularities, and $V'$ is smooth, and
    \item there exist two $\Rr$-divisors $B_Y$ and $E$ on $Y$, such that
    \begin{itemize}
    \item $K_Y+B_Y+\Mm_Y=f^*(K_X+B+\Mm_X)+E$,
    \item $B_Y\geq 0$, $E\geq 0$, and $B_Y\wedge E=0$,
    \item $(Y,B_Y)$ is lc quasi-smooth, and any lc center of $(Y,B_Y,\Mm)$ on $X$ is an lc center of $(X,B,\Mm)$.
    \end{itemize}
\end{itemize}
In particular, $(Y,B_Y,\Mm)$ is $\Qq$-factorial NQC lc and $Y$ is klt. Since $\varphi$ is birational, by Lemma \ref{lem: property of numerical and Iitaka dimension}(3) we obtain
$$\kappa_{\sigma}(Y/V',K_Y+B_Y+\Mm_Y)=\kappa_{\sigma}(Y/V,K_Y+B_Y+\Mm_Y)=\kappa_{\sigma}(X/V,K_X+B+\Mm_X)=0$$
and
$$\kappa_{\iota}(Y/V',K_Y+B_Y+\Mm_Y)=\kappa_{\iota}(Y/V,K_Y+B_Y+\Mm_Y)=\kappa_{\iota}(X/V,K_X+B+\Mm_X)=0.$$
By Lemma \ref{lem: has19 3.2 step 3 abu ver}, we may run a $(K_Y+B_Y+\Mm_Y)$-MMP$/V'$ with scaling of a general ample$/V'$ divisor $A$ on $Y$, which terminates with a good minimal model $(X',B',\Mm)/V'$ of $(Y,B_Y,\Mm)/V'$ such that $K_{X'}+B'+\Mm_{X'}\sim_{\Rr,V'}0$. Let $\pi': X'\rightarrow V'$ be the induced contraction. 

We show that $(X',B',\Mm)/U,\pi',\varphi$ satisfy our requirements. (1)(2) follow from our construction. 

Let $p: W'\rightarrow Y$ and $q: W'\rightarrow X'$ be a resolution of indeterminacies of the induced map $Y\dashrightarrow X'$ such that $p$ is a log resolution of $(Y,B_Y)$.
\begin{center}$\xymatrix{
& W'\ar@{->}[d]^{q}\ar@{->}[dl]_{p} & \\
X'\ar@{->}[dr]_{\pi'}& Y\ar@{-->}[l]\ar@{->}[r]^{f}\ar@{->}[d]_{\pi_Y} & X\ar@{->}[d]^{\pi}\\
& V'\ar@{->}[r]^{\varphi} & V
}$
\end{center}

Then we have
$$p^*(K_Y+B_Y+\Mm_Y)=q^*(K_{X'}+B'+\Mm_{X'})+F$$
for some $F\geq 0$ that is exceptional over $X'$. Let $B_{W'}:=p^{-1}_*B_Y+\Exc(p)$, then $(W',B_{W'},\Mm)$ is a log smooth model of $(Y,B_Y,\Mm)$ and $(X',B',\Mm)$.

Since $K_Y+B_Y+\Mm_Y=f^*(K_{X}+B+\Mm_{X})+E$, by \cite[Theorem 3.14]{HL21a}, $(X,B,\Mm)/U$ has a good minimal model if and only if $(Y,B_Y,\Mm)/U$ has a good minimal model, if and only if $(W',B_{W'},\Mm)/U$ has a good minimal model, if and only if $(X',B',\Mm)/U$ has a good minimal model, hence (3).

By \cite[Lemmas 3.10, 3.17]{HL21a}, a g-pair $(X'',B'',\Mm)/U$ is a weak lc model of $(X,B,\Mm)/U$ if and only if $(X'',B'',\Mm)/U$ is a weak lc model of $(W',B_{W'},\Mm)/U$, if and only if $(X'',B'',\Mm)/U$ is a weak lc model of $(X',B',\Mm)/U$, hence (4).

Let $D$ be an lc place of $(X',B',\Mm)$. Since $Y\dashrightarrow X'$ is a $(K_{Y}+B_Y+\Mm_Y)$-MMP$/V'$, $D$ is an lc place of $(Y,B_Y,\Mm)$, hence an lc place of $(X,B,\Mm)$. Thus if all lc centers of $(X,B,\Mm)$ dominate $V$, then all lc centers of $(X',B',\Mm)$ dominate $V$, hence all lc centers of $(X',B',\Mm)$ dominate $V'$ as $\varphi$ is birational, and we have (5). 

Finally, by Lemma \ref{lem: property of numerical and Iitaka dimension}(3) we obtain
\begin{align*}
    \kappa_{\sigma}(X/U,K_X+B+\Mm_X)&=\kappa_{\sigma}(Y/U,K_Y+B_Y+\Mm_Y)=\kappa_{\sigma}(W'/U,p^*(K_Y+B_Y+\Mm_Y))\\
    &=\kappa_{\sigma}(W'/U,q^*(K_{X'}+B'+\Mm_{X'})+F)\\
    &=\kappa_{\sigma}(X'/U,K_{X'}+B'+\Mm_{X'})
\end{align*}
and
\begin{align*}
    \kappa_{\iota}(X/U,K_X+B+\Mm_X)&=\kappa_{\iota}(Y/U,K_Y+B_Y+\Mm_Y)=\kappa_{\iota}(W'/U,p^*(K_Y+B_Y+\Mm_Y))\\
    &=\kappa_{\iota}(W'/U,q^*(K_{X'}+B'+\Mm_{X'})+F)\\
    &=\kappa_{\iota}(X'/U,K_{X'}+B'+\Mm_{X'}),
\end{align*}
and we get (6).
\end{proof}

\begin{prop}[cf.{ \cite[Lemma 3.10]{Has22a}}]\label{prop: prop 3.4 has19 abu ver}
Let $(X,B,\Mm)/U$ be an NQC lc g-pair and $\pi: X\rightarrow V$ a contraction over $U$, such that
\begin{itemize}
    \item $V$ is normal quasi-projective,
    \item $\kappa_{\sigma}(X/V,K_X+B+\Mm_X)=\kappa_{\iota}(X/V,K_X+B+\Mm_X)=0$ and $\kappa_{\sigma}(X/U,K_X+B+\Mm_X)=\dim V-\dim U$, and 
    \item all lc centers of $(X,B,\Mm)$ dominate $V$.
\end{itemize}
Then:
\begin{enumerate}
    \item $(X,B,\Mm)/U$ has a good minimal model, and
    \item Let $(\bar X,\bar B,\Mm)/U$ be a good minimal model of $(X,B,\Mm)/U$ and $\bar X\rightarrow\bar V$ is the contraction over $U$ induced by $K_{\bar X}+\bar B+\Mm_{\bar X}$. Then all lc centers of $(\bar X,\bar B,\Mm)$ dominate $\bar V$.
\end{enumerate}
\end{prop}
\begin{proof}
By Theorem \ref{thm: HL21a v2 4.1 abu ver}, there exists a $\Qq$-factorial NQC dlt g-pair $(X',B',\Mm)/U$, a contraction $\pi': X'\rightarrow V'$ over $U$, and a birational projective morphism $\varphi: V'\rightarrow V$ over $U$ such that 
\begin{itemize}
    \item $X'$ is birational to $X$ and $V'$ is smooth,
    \item $K_{X'}+B'+\Mm_{X'}\sim_{\mathbb R,V'}0$. In particular, $\kappa_{\sigma}(X'/V',K_{X'}+B'+\Mm_{X'})=0$ by Lemma \ref{lem: property of numerical and Iitaka dimension}(5),
    \item $(X,B,\Mm)/U$ has a good minimal model if and only if $(X',B',\Mm)/U$ has a good minimal model,
    \item any weak lc model of $(X,B,\Mm)/U$ is a weak lc model of $(X',B',\Mm)/U$, and any weak lc model of $(X',B',\Mm)/U$ is a  weak lc model of $(X,B,\Mm)/U$,
    \item all lc centers of $(X',B',\Mm)$ dominate $V'$, and
    \item $\kappa_{\sigma}(X'/U,K_{X'}+B'+\Mm_{X'})=\kappa_{\sigma}(X/U,K_X+B+\Mm_X)=\dim V-\dim U=\dim V'-\dim U$.
\end{itemize}
\begin{center}$\xymatrix{
X'\ar@{->}[d]_{\pi'}\ar@{-->}[rr]& & X\ar@{->}[d]^{\pi}\\
V'\ar@{->}[rr]^{\varphi}\ar@{->}[dr] & & V\ar@{->}[dl]\\
& U &
}$
\end{center}
\begin{claim}\label{claim: check last part of prop 3.4 has19 for g pair}
Assume that $(X',B',\Mm)/U$ has a good minimal model $(\bar X',\bar B',\Mm)/U$, $\bar X'\rightarrow \bar V'$ is the contraction over $U$ induced by $K_{\bar X'}+\bar B'+\Mm_{\bar X'}$, and all lc centers of $(\bar X',\bar B',\Mm)$ dominate $\bar V'$. Then Proposition \ref{prop: prop 3.4 has19 abu ver}(2) holds for $(X,B,\Mm)/U$.
\end{claim}
\begin{proof}
Let $(\bar X,\bar B,\Mm)/U$ be a good minimal model of $(X,B,\Mm)/U$. Then $(\bar X,\bar B,\Mm)/U$ is a weak lc model of $(X',B',\Mm)/U$. Since $(\bar X',\bar B',\Mm)/U$ is also a weak lc model of  $(X',B',\Mm)/U$, by \cite[Lemma 3.9(1)]{HL21a}, we may take a resolution of indeterminacies $p: W\rightarrow\bar X$ and $q: W\rightarrow\bar X'$ of the induced birational map $\bar X\dashrightarrow\bar X'$ such that 
$$p^*(K_{\bar X}+\bar B+\Mm_{\bar X})=q^*(K_{\bar X'}+\bar B'+\Mm_{\bar X'}).$$
Then:
\begin{itemize}
    \item $K_{\bar X}+\bar B+\Mm_{\bar X}$ is semi-ample$/U$, and if we let $\bar X\rightarrow\bar V$ be the contraction over $U$ induced by $K_{\bar X}+\bar B+\Mm_{\bar X}$, then $\bar V=\bar V'$ since they are defined by the same linear series. 
    \item Any lc center of $(\bar X,\bar B,\Mm)$ is an lc center of $(\bar X',\bar B',\Mm)$, and any lc center of $(\bar X',\bar B',\Mm)$ is an lc center of $(\bar X,\bar B,\Mm)$. In particular, since all lc centers of $(\bar X',\bar B',\Mm)$ dominate $\bar V'=\bar V$, all lc centers of  $(\bar X,\bar B,\Mm)$ dominate $\bar V$.
\end{itemize}
The claim is proved.
\end{proof}
\noindent\textit{Proof of Proposition \ref{prop: prop 3.4 has19 abu ver} continued}. By Claim \ref{claim: check last part of prop 3.4 has19 for g pair}, we may replace $(X,B,\Mm),V$ and $\pi$ with $(X',B',\Mm),V'$ and $\pi'$ respectively, and assume that $V$ is smooth and $K_X+B+\Mm_X\sim_{\mathbb R,V}0$. By Theorem \ref{thm: numerical generalized canonical bundle formula}, there exists an NQC klt g-pair $(V,B_V,\Mm^V)/U$ such that 
$$K_X+B+\Mm_X\sim_{\mathbb R}\pi^*(K_V+B_V+\Mm^V_V).$$
By Lemma \ref{lem: property of numerical and Iitaka dimension}(4)(5), we have
$$\kappa_{\sigma}(V/U,K_V+B_V+\Mm^V_V)=\kappa_{\sigma}(X/U,K_X+B+\Mm_X)=\dim V-\dim U.$$
By Lemma \ref{lem: property of numerical and Iitaka dimension}(1), $K_V+B_V+\Mm^V_V$ is big$/U$. By \cite[Lemma 4.4(2)]{BZ16}, we may run a $(K_V+B_V+\Mm^V_V)$-MMP$/U$ with scaling of some general ample$/U$ divisor $A$, which terminates with a good minimal model $(\widehat V,B_{\widehat V},\Mm^V)/U$ of $(V,B_V,\Mm^V)/U$. Let $\phi: V\dashrightarrow\widehat V$ be the induced morphism, and let $g: \tilde V\rightarrow V$ and $\widehat g: \tilde V\rightarrow\widehat V$ be a common resolution such that $\widehat g=\phi\circ g$. Then
$$g^*(K_V+B_V+\Mm^V_V)=\widehat g^*(K_{\widehat V}+B_{\widehat V}+\Mm^V_{\widehat V})+F.$$
for some $\widehat g$-exceptional $\Rr$-divisor $F\ge0$ on $\tilde V$. Let $h: W\rightarrow X$ be a log resolution of $(X,\Supp B)$ such that $\Mm$ descends to $W$ and the induced map $\pi_W: W\rightarrow\tilde V$ is a morphism. 

\begin{center}$\xymatrix{
Z\ar@/^2pc/[rr]^{\widehat{f}}\ar@{->}[d]_{f} & W\ar@{->}[d]_{h}\ar@{->}[dr]^{\pi_W}\ar@{-->}[r] & \widehat{W}\ar@/^2pc/[dd]^{\pi_{\widehat W}}\\
\bar X\ar@{->}[d]& X\ar@{->}[d]_{\pi}\ar@{-->}[l]& \tilde V\ar@{->}[dl]_{g}\ar@{->}[d]^{\widehat g} \\
\bar V& V\ar@{-->}[r]^{\phi}\ar@{-->}[l] & \widehat V\ar@/^2pc/[ll]
}$
\end{center}

\vspace{5mm}

By \cite[Lemma 3.6]{HL21a}, there exists a proper log smooth model $(W,B_W,\Mm)$ of $(X,B,\Mm)$. In particular,
$$K_W+B_W+\Mm_W=h^*(K_X+B+\Mm_X)+E$$
for some $h$-exceptional $\Rr$-divisor $E\geq 0$. 
Then
\begin{align*}
    K_W+B_W+\Mm_W&=h^*(K_X+B+\Mm_X)+E\sim_{\mathbb R}(\pi\circ h)^*(K_V+B_V+\Mm^V_V)+E\\
    &=\pi_W^*g^*(K_V+B_V+\Mm^V_V)+E=\pi_W^*\widehat g^*(K_{\widehat V}+B_{\widehat V}+\Mm^V_{\widehat V})+\pi_W^*F+E.
\end{align*}
Since $E$ is exceptional over $X$, $E$ is very exceptional over $V$ (see \cite[Paragraph after Definition 3.1]{Bir12}). Since $\phi$ is a birational contraction, $E$ is very exceptional over $\widehat V$. Since $F$ is exceptional over $\widehat V$, $\pi_W^*F$ is very exceptional over $\widehat V$. Therefore 
$$K_W+B_W+\Mm_W\sim_{\mathbb R,\widehat V}\pi_W^*F+E$$
is very exceptional over $\widehat V$. By Lemma \ref{lem: rlinear version of HL22 3.8}, we may run a $(K_W+B_W+\Mm_W)$-MMP$/\widehat V$ with scaling of a general ample$/\widehat V$ divisor which terminates with a good minimal model $(\widehat W,B_{\widehat W},\Mm)/\widehat V$ such that $K_{\widehat W}+B_{\widehat W}+\Mm_{\widehat W}\sim_{\Rr,\widehat V}0$ and the induced birational map $W\dashrightarrow\widehat W$ exactly contracts $\Supp(\pi_W^*F+E)$. In particular, let $\pi_{\widehat W}: \widehat W\rightarrow\widehat  V$ be the induced morphism, then $$K_{\widehat W}+B_{\widehat W}+\Mm_{\widehat W}\sim_{\Rr}\pi_{\widehat W}^*(K_{\widehat V}+B_{\widehat V}+\Mm^V_{\widehat V}).$$
Since $(\widehat V,B_{\widehat V},\Mm^V)/U$ is a good minimal model of $(V,B_V,\Mm^V)/U$, $K_{\widehat V}+B_{\widehat V}+\Mm^V_{\widehat V}$ is semi-ample$/U$, hence $K_{\widehat W}+B_{\widehat W}+\Mm_{\widehat W}$ is semi-ample$/U$. Thus $(\widehat W,B_{\widehat W},\Mm)/U$ is a good minimal model of $(W,B_W,\Mm)/U$. By \cite[Lemma 3.10]{HL21a}, $(\widehat W,B_{\widehat W},\Mm)/U$ is a good minimal model of $(X,B,\Mm)/U$, which implies (1).

Let $(\bar X,\bar B,\Mm)/U$ be a good minimal model of $(X,B,\Mm)/U$. By \cite[Lemma 3.9(1)]{HL21a}, there exists a resolution $f: Z\rightarrow \bar X$ and $\widehat f: Z\rightarrow\widehat W$ of indeterminacies of the induced birational map $\bar X\dashrightarrow\widehat W$ such that 
$$f^*(K_{\bar X}+\bar B+\Mm_{\bar X})=\widehat f^*(K_{\widehat W}+B_{\widehat W}+\Mm_{\widehat W}).$$ In particular, any lc place of $(\bar X,\bar B,\Mm)$ is an lc place of $(\widehat W,B_{\widehat W},\Mm)$, hence an lc place of $(W,B_W,\Mm)$, and thus an lc place of $(X,B,\Mm)$ by \cite[Lemma 3.7]{HL21a}. Therefore, any lc place of $(\bar X,\bar B,\Mm)$ dominates $V$. Moreover, since 
$$
f^*(K_{\bar X}+\bar B+\Mm_{\bar X})\sim_{\Rr}\widehat f^*\circ\pi_{\widehat W}^*(K_{\widehat V}+B_{\widehat V}+\Mm^V_{\widehat V}),
$$
the contraction $ Z\rightarrow\bar V$ induced by  $f^*(K_{\bar X}+\bar B+\Mm_{\bar X})$ factors through $\widehat V$, and the induced morphism $\widehat V\rightarrow\bar V$ is actually given by the big$/U$ semi-ample$/U$ $\Rr$-divisor $K_{\widehat V}+B_{\widehat V}+\Mm^V_{\widehat V}$.
In particular, the induced map $V\dashrightarrow\bar V$ is birational. Thus all lc places of $(\bar X,\bar B,\Mm)$ dominate $\bar V$, hence all lc centers of  $(\bar X,\bar B,\Mm)$ dominate $\bar V$, which implies (2).
\end{proof}

\section{Applications of the Nakayama-Zariski decomposition}

This section is similar to \cite[Section 3, before Theorem 3.14]{Has22a}.

\begin{lem}[{cf. \cite[Lemma 3.5]{Has22a}}]\label{lem: Has22a 3.5 rel} 
Let $(X,B,\Mm)/U$ and $(X',B',\Mm)/U$ be NQC dlt g-pairs with a birational map $\phi: X\dashrightarrow X'$ over $U$ such that $\phi_*\Mm=\Mm$. Let $S$ and $S'$ be lc centers of $(X,B,\Mm)$ and $(X',B',\Mm)$ respectively, such that $\phi$ is an isomorphism near the generic point of $S$, and $\phi|_S: S\dashrightarrow S'$ defines a birational map$/U$. Suppose that
\begin{enumerate}
    \item $K_{X}+B+\Mm_X$ is pseudo-effective$/U$,
    \item for any prime divisor $D'$ on $X'$, $a(D',X',B',\Mm)\leq a(D',X,B,\Mm)$, and
    \item for every prime divisor $P$ over $X$ such that $a(P,X,B,\Mm)<1$ and $\Center_{X}(P)\cap S\not=\emptyset$, then $\sigma_{P}(X/U,K_X+B+\Mm_X)=0$.
\end{enumerate}
Let $(S,B_S,\Mm^S)/U$ and $(S',B_{S'},\Mm^S)/U$ be the dlt g-pairs induced by adjunction 
$$
K_S+B_S+\Mm^S_S:=(K_X+B+\Mm_X)|_S~,~K_{S'}+B_{S'}+\Mm^S_{S'}:=(K_{X'}+B'+\Mm_{X'})|_{S'}.
$$
Then 
$$
a(Q,S',B_{S'},\Mm^S)\leq a(Q,S,B_S,\Mm^S)
$$ for all prime divisors $Q$ on $S'$. 
\end{lem} 
\begin{proof}
The proof follows exactly the same lines as \cite[Proof of Lemma 3.5]{Has22a} except the following two places:
\begin{itemize}
    \item \cite[Page 13, Line 30]{Has22a} cites \cite[Remark 2.3(1)]{Has20}. We shall replace \cite[Remark 2.3(1)]{Has20} with Lemma \ref{lem: nz basic properties}(1).
    \item \cite[Page 13, Line 32]{Has22a} cites \cite[Remark 2.3(3)]{Has20}. We shall replace \cite[Remark 2.3(3)]{Has20} with Lemma \ref{lem: nz keep under pullback}(2).
\end{itemize}
Therefore, we shall omit the details of the proof to avoid redundance.
\end{proof}

\begin{lem}[{cf. \cite[Lemma 5.3]{HMX18}}]\label{lem: Has22a 3.6 rel dlt} Let $(X,B_1,\Mm)/U$ and $(X,B_2,\Mm)/U$ be $\Qq$-factorial NQC dlt g-pairs such that $K_X+B_1+\Mm_X$ is pseudo-effective$/U$ and 
$$
0\leq B_1-B_2\leq N_{\sigma}(X/U,K_X+B_1+\Mm_X).
$$ Then $(X,B_1,\Mm)/U$ has a log minimal model (resp. good minimal model) if and only if $(X,B_2,\Mm)/U$ has a log minimal model (resp. good minimal model).  
\end{lem}
\begin{proof}

First we assume that $(X,B_1,\Mm)/U$ has a log minimal model (resp. good minimal model). By \cite[Theorem 2.24]{HL21a}, we may run a $(K_X+B_1+\Mm_X)$-MMP$/U$ which terminates with a log minimal model (resp. good minimal model) $(X',B',\Mm)/U$ with induced birational map $\phi: X\dashrightarrow X'$ over $U$. By Lemmas \ref{lem: nz pre properties}(1) and \ref{lem: nz for lc divisor}(2), $\phi$ contracts every component of $\Supp N_{\sigma}(X/U,K_X+B_1+\Mm_X)$. Thus $B'$ is also the strict transform of $B_2$ on $X'$.

Let $p: W\rightarrow X$ and $q: W\rightarrow X'$ be a resolution of indeterminacies of $\phi$, and write
$$p^*(K_X+B_1+\Mm_X)=q^*(K_{X'}+B'+\Mm_{X'})+E$$
for some effective $q$-exceptional $\Rr$-divisor $E$ on $W$.
Then by Lemmas \ref{lem: nz pre properties}(1) and \ref{lem: nz keep under pullback}(2)(3), $N_{\sigma}(X/U,K_X+B_1+\Mm_X)=p_*E$ is well-defined as a divisor. Let $F:=E-p^*(B_1-B_2)$. Then
$$F\geq E-p^*N_{\sigma}(X/U,K_X+B_1+\Mm_X)=E-p^*p_*E$$
and
$$p^*(K_X+B_2+\Mm_X)=q^*(K_{X'}+B'+\Mm_{X'})+F.$$
Since $E-p^*p_*E$ is $p$-exceptional, $p_*F\geq 0$. By the negativity lemma, $F\geq 0$. Thus $(X',B',\Mm)/U$ is a weak lc model of $(X,B_2,\Mm)/U$. By \cite[Lemmas 3.9(2), 3.15]{HL21a}, $(X,B_2,\Mm)/U$ has a log minimal model (resp. good minimal model).

Now we assume that $(X,B_2,\Mm)/U$ has a log minimal model (resp. good minimal model). By \cite[Theorem 2.24]{HL21a}, we may run a $(K_X+B_2+\Mm_X)$-MMP$/U$ which terminates with a log minimal model (resp. good minimal model) $(X',B',\Mm)/U$ with induced birational map $\phi: X\dashrightarrow X'$ over $U$. Let $C:=B_1-B_2$. Then $\phi$ is also a  $(K_X+B_2+\epsilon C+\Mm_X)$-MMP$/U$ for any $0<\epsilon\ll 1$ by Lemma \ref{lem: still an mmp under perturbation}. Let $C'$ be the strict transform of $C$  on $X'$. By Lemma \ref{lem: nz basic properties}(5), we obtain
$$N_{\sigma}(X/U,K_X+B_2+\epsilon C+\Mm_X)+(1-\epsilon)C=N_{\sigma}(X/U,K_X+B_1+\Mm_X),$$
and
$$N_{\sigma}(X/U,K_X+B_2+\Mm_X)+C=N_{\sigma}(X/U,K_X+B_1+\Mm_X)$$
for any $\epsilon\in[0,1]$. Therefore, we have
$$
N_{\sigma}(X/U,K_X+B_2+\epsilon C+\Mm_X)=N_{\sigma}(X/U,K_X+B_2+\Mm_X)+\epsilon C
$$ 
for any $\epsilon\in[0,1]$. Hence, if $\epsilon\in(0,1]$, then
$$
\Supp N_{\sigma}(X/U,K_X+B_2+\epsilon C+\Mm_X)=\Supp N_{\sigma}(X/U,K_X+B_1+\Mm_X),
$$
since they are both equal to $\Supp N_{\sigma}(X/U,K_X+B_2+\Mm_X)\cup\Supp C$. Moreover, by \cite[Lemma 3.21]{HL22}, we may pick $0<\epsilon\ll 1$ such that any partial $(K_{X'}+B_2'+\epsilon C'+\Mm_{X'})$-MMP$/U$ is $(K_{X'}+B'+\Mm_{X'})$-trivial$/U$. We run a $(K_{X'}+B'+\epsilon C'+\Mm_{X'})$-MMP$/U$ with scaling of an ample$/U$ $\Rr$-divisor. By Lemma \ref{lem: limit movable r divisors gpairs}, after finitely many steps we get a birational map $\psi: X'\dashrightarrow X''$ such that $K_{X''}+B''+\epsilon C''+\Mm_{X''}$ is a movable$/U$ $\Rr$-divisor, where $B''$ and $C''$ are the strict transforms of $B'$ and $C'$ on $X''$ respectively. By Lemma \ref{lem: nz for lc divisor}(2), the set of $(\psi\circ\phi)$-exceptional divisors is exactly $\Supp N_{\sigma}(X/U,K_X+B_2+\epsilon C+\Mm_X)=\Supp N_{\sigma}(X/U,K_X+B_1+\Mm_X)$. Thus, by assumption, $\Supp C=\Supp (B_1-B_2)$ is also $(\psi\circ\phi)$-exceptional. Then $C''=0$ since it is the pushforward of $C$ to $X''$, hence $B''$ is also the strict transform of $B_1$ on $X''$ and $K_{X''}+B''+\Mm_{X''}$ is nef$/U$ (resp. semi-ample$/U$) by construction. By Lemma \ref{lem: hmx18 2.7.3 gpair rel}, $(X'',B'',\Mm)/U$ is a log minimal model of $(X,B_1,\Mm)/U$. The lemma follows from \cite[Lemma 3.9(2)]{HL21a}.
\end{proof}

\begin{lem}[{cf. \cite[Lemma 3.6]{Has22a}}]\label{lem: Has22a 3.6 rel} Let $(X,B,\Mm)/U$ and $(Y,B_Y,\Mm)/U$ be NQC lc g-pairs and $f: Y\rightarrow X$ a projective birational morphism such that
\begin{enumerate}
    \item $K_X+B+\Mm_X$ is pseudo-effective$/U$, and
    \item for any prime divisor $D$ on $Y$,
    $$0\leq a(D,Y,B_Y,\Mm)-a(D,X,B,\Mm)\leq\sigma_D(X/U,K_X+B+\Mm_X).$$
\end{enumerate}
Then $K_Y+B_Y+\Mm_Y$ is pseudo-effective$/U$. Moreover, $(X,B,\Mm)/U$ has a log minimal model (resp. good minimal model) if and only if $(Y,B_Y,\Mm)/U$ has a log minimal model (resp. good minimal model).
\end{lem}
\begin{proof}
The assumptions imply that 
$$
0\le f^*(K_X+B+\Mm_X)-(K_Y+B_Y+\Mm_Y)\le N_\sigma(Y/U,f^*(K_X+B+\Mm_X))
$$
and then $K_Y+B_Y+\Mm_Y$ is pseudo-effective$/U$ by Lemma \ref{lem: nz keep under pullback}(4).

Let $g: W\rightarrow Y$ be a log resolution of $(Y,\Supp B_Y)$ such that $\Mm$ descends to $W$ and $h:=f\circ g$ is a log resolution of $(X,\Supp B)$. Let $B_W:=h^{-1}_*B+\Supp\Exc(h)$ and $B'_W:=g^{-1}_*B_Y+\Supp\Exc(g)$. Then we have
$$K_W+B_W+\Mm_W=h^*(K_X+B+\Mm_X)+E$$
for some $E_W\geq 0$ that is exceptional$/X$. By Lemma \ref{lem: nz keep under pullback}(1)(2), 
$$\sigma_P(W/U,K_W+B_W+\Mm_W)=\sigma_P(X/U,K_X+B+\Mm_X)+\mult_PE$$
for any prime divisor $P$ on $W$.
\begin{claim}\label{claim: has 3.6 change to dlt}
For any prime divisor $P$ on $W$, we have
$$0\leq a(P,W,B_W',\Mm)-a(P,W,B_W,\Mm)\leq\sigma_P(W/U,K_W+B_W+\Mm_W).$$
\end{claim}
Grant Claim \ref{claim: has 3.6 change to dlt} for the time being. By Claim \ref{claim: has 3.6 change to dlt} and \cite[Theorem 3.14]{HL21a}, possibly replacing $(X,B,\Mm)/U$ and $(Y,B_Y,\Mm)/U$ with $(W,B_W,\Mm)/U$ and $(W,B_W',\Mm)/U$ respectively, we may assume that $(X,B,\Mm)$ and $(Y,B_Y,\Mm)$ are $\Qq$-factorial dlt and $X=Y$. The lemma follows from Lemma \ref{lem: Has22a 3.6 rel dlt}.
\end{proof}

\begin{proof}[Proof of Claim \ref{claim: has 3.6 change to dlt}]
For any prime divisor $P$ on $W$, one of the following cases holds:

\medskip

\noindent\textbf{Case 1}. $P$ is not exceptional over $X$. In this case 
$$a(P,W,B_W',\Mm)-a(P,W,B_W,\Mm)=a(P,Y,B_Y,\Mm)-a(P,X,B,\Mm)$$
and the claim follows.

\medskip

\noindent\textbf{Case 2}. $P$ is exceptional over $X$ but not exceptional over $Y$. In this case $a(P,W,B_W,\Mm)=0$, $a(P,W,B_W',\Mm)=a(P,Y,B_Y,\Mm)$, and $a(P,X,B,\Mm)=\mult_PE$, so
$$0\leq a(P,Y,B_Y,\Mm)=a(P,W,B_W',\Mm)-a(P,W,B_W,\Mm),$$
and
\begin{align*}
    a(P,Y,B_Y,\Mm)&\leq\sigma_P(X/U,K_X+B+\Mm_X)+a(P,X,B,\Mm)\\
    &=\sigma_P(X/U,K_X+B+\Mm_X)+\mult_PE=\sigma_P(W/U,K_W+B_W+\Mm_W)
\end{align*}
and the claim follows.

\medskip

\noindent\textbf{Case 3}. $P$ is exceptional over $Y$. In this case $a(P,W,B_W,\Mm)=a(P,W,B'_W,\Mm)=0$, and the claim follows.
\end{proof}

\begin{lem}[{cf. \cite[Lemma 3.8]{Has22a}}]\label{lem: Has22a 3.8 rel}
Let $(X,B,\Mm)/U$ be an NQC lc g-pair with induced morphism $\pi: X\rightarrow U$ such that $U$ is quasi-projective. Let $S$ be a subvariety of $X$, and
$$(X,B,\Mm):=(X_0,B_0,\Mm)\dashrightarrow (X_1,B_1,\Mm)\dashrightarrow\dots\dashrightarrow (X_n,B_n,\Mm)\dashrightarrow\dots$$
a $(K_X+B+\Mm_X)$-MMP$/U$ with scaling of an 
$\Rr$-divisor $A\geq 0$. Let
$$\lambda_i:=\inf\{t\geq 0\mid K_{X_i}+B_i+\Mm_{X_i}+tA_i\text{ is nef}/U\}$$
be the scaling numbers, where $A_i$ is the strict transform of $A$ on $X_i$. Suppose that
\begin{itemize}
    \item each step of this MMP is an isomorphism on a neighborhood of $S$, and
    \item $\lim_{i\rightarrow+\infty}\lambda_i=0$.
\end{itemize}
Then 
\begin{enumerate}
    \item for any $\pi$-ample $\Rr$-divisor $H$ on $X$ and any closed point $x\in S$, there exists an $\Rr$-divisor $E$ such that $0\leq E\sim_{\mathbb{R},U}K_X+B+\Mm_X+H$ and $x\not\in\Supp E$, and
    \item  for any prime divisor $P$ over $X$ such that $\Center_XP\cap S\not=\emptyset$, $\sigma_P(X/U,K_X+B+\Mm_X)=0$.
\end{enumerate} 
\end{lem}
\begin{proof}
(1) follows from \cite[Lemma 3.8]{Has22a} and (2) follows from (1) and Lemma \ref{lem: nz basic properties}(4).
\end{proof}

\begin{lem}[{cf. \cite[Lemma 3.9]{Has22a}}]\label{lem: Has22a 3.9 rel} Let $(X,B,\Mm)/U$ and $(X',B',\Mm)/U$ be two NQC lc g-pairs and $\phi: X\dashrightarrow X'$ a birational map such that $\phi_*\Mm=\Mm$. Suppose that
\begin{itemize}
    \item $a(P,X,B,\Mm)\leq a(P,X',B',\Mm)$ for any prime divisor $P$ on $X$, and
    \item $a(P',X',B',\Mm)\leq a(P',X,B,\Mm)$ for any prime divisor $P'$ on $X'$. 
\end{itemize}
Then
\begin{enumerate}
    \item $K_X+B+\Mm_{X}$ is abundant$/U$ if and only if $K_{X'}+B'+\Mm_{X'}$ is abundant$/U$, and
    \item $(X,B,\Mm)/U$ has a log minimal model (resp. good minimal model) if and only if $(X',B',\Mm)/U$ has a log minimal model (resp. good minimal model).
\end{enumerate}
\end{lem}

\begin{proof}
Let $p: W\rightarrow X$ and $q: W\rightarrow X'$ be a resolution of indeterminacies such that $\Mm$ descends to $W$, $p$ is a log resolution of $(X,\Supp B)$, and $q$ is a log resolution of $(X',\Supp B')$. Let
$$B_W:=\sum_{D\text{ is a prime divisor on }W}\max\{1-a(D,X,B,\Mm),1-a(D,X',B',\Mm),0\}D.$$
Then $(W,B_W,\Mm)$ is lc and $(W,B_W)$ is log smooth. By construction, there exists a $p$-exceptional $\Rr$-divisor $E\geq 0$ and a $q$-exceptional $\Rr$-divisor $F\geq 0$ such that
$$E+p^*(K_X+B+\Mm_X)=K_W+B_W+\Mm_W=q^*(K_{X'}+B'+\Mm_{X'})+F.$$
(1) follows from Lemma \ref{lem: property of numerical and Iitaka dimension}(3) and (2) follows from \cite[Theorem 3.14]{HL21a}.
\end{proof}

\section{A special log minimal model}

The purpose of this section is to prove Theorem \ref{thm: Has22a 3.14 rel ver} and Theorem \ref{thm: Has22a 3.15 rel ver}, which are analogues of \cite[Theorem 3.14 and Theorem 3.15]{Has22a} in the relative setting. 

\begin{thm}[{cf. \cite[Theorem 3.14]{Has22a}}]\label{thm: Has22a 3.14 rel ver}
Let $(X,B,\Mm)/U$ be an NQC dlt g-pair such that
\begin{itemize}
    \item $K_X+B+\Mm_X$ is pseudo-effective$/U$ and abundant$/U$,
    \item for any lc center $S$ of $(X,B,\Mm)$, $(K_X+B+\Mm_X)|_S$ is nef$/U$, and
    \item for any prime divisor $P$ over $X$ such that $a(P,X,B,\Mm)<1$ and $\Center_XP\cap\Nklt(X,B,\Mm)\not=\emptyset$, $\sigma_P(X/U,K_X+B+\Mm_X)=0$.
\end{itemize}
Then $(X,B,\Mm)/U$ has a log minimal model.
\end{thm}

\begin{proof}
We divide the proof in six steps.
\medskip

\noindent\textbf{Step 1}. In this step we show that we may replace $(X,B,\Mm)$ with a $\Qq$-factorial dlt model and find two $\Rr$-divisors $G\geq 0,H\geq 0$, and a real number $1>t_0>0$ such that
\begin{itemize}
    \item[(I)] $K_X+B+\Mm_X\sim_{\Rr,U}G+H$,
    \item[(II)] $\Supp G\subset\Supp\lfloor B\rfloor$, and
    \item[(III)] for any $t\in (0,t_0]$, the following hold:
    \begin{enumerate}
        \item[(III.1)] $(X,B+tH,\Mm)/U$ is dlt, $N_{\sigma}(X/U,K_X+B+tH+\Mm_X)$ is well-defined as a divisor and
        $\Supp N_{\sigma}(X/U,K_X+B+tH+\Mm_X)$ 
        does not depend on $t$, and
        \item[(III.2)] $(X,B-tG,\Mm)/U$ has a good minimal model.
    \end{enumerate}
\end{itemize}

Since $K_X+B+\Mm_X$ is pseudo-effective$/U$ and abundant$/U$,  $K_{X}+B+\Mm_{X}\sim_{\Rr,U}D\geq 0$ for some $\Rr$-divisor $D$ on $X$. Let $X\dashrightarrow V$ be the Iitaka fibration$/U$ associated to $D$. Then $\dim V-\dim U=\kappa_{\sigma}(X/U,K_X+B+\Mm_X)$. Let $h: W\rightarrow X$ be a log resolution of $(X,\Supp B)$ such that $\Mm$ descends to $W$ and the induced map $\psi: W\dashrightarrow V$ is a morphism. Then we may write
$$K_W+B_W+\Mm_W=h^*(K_X+B+\Mm_X)+E$$
such that $B_W\geq 0,E\geq 0$, and $B_W\wedge E=0$. Notice that $(W,B_W,\Mm)$ is a log smooth model of $(X,B,\Mm)$. By Lemma \ref{lem: iitaka fibration numerical abundant divisor gpair dimension}, 
\begin{itemize}
    \item[(i)]  $\kappa_{\sigma}(W/U,K_{W}+B_W+\Mm_W)=\dim V-\dim U$ and $\kappa_{\sigma}(W/V,K_{W}+B_W+\Mm_{W})=0$. 
\end{itemize}
Thus by construction $K_{W}+B_W+\Mm_{W}$ is $\Rr$-linearly equivalent$/U$ to the sum of an effective $\Rr$-divisor and the pullback of an ample$/U$ $\Rr$-divisor on $V$. In particular, we may find $0\leq D_W\sim_{\Rr,U}K_W+B_W+\Mm_W$ such that $\Supp D_W$ contains all lc centers of $(W,B_W,\Mm)$ that are vertical over $V$. 

Let $(\bar X,\bar B,\Mm)$ be a proper log smooth model of $(W,B_W,\Mm)$ with induced morphism $g: \bar X\rightarrow W$ such that $g$ is a log resolution of $(W,B_W+D_W)$, and
$$K_{\bar X}+\bar B+\Mm_{\bar X}=g^*(K_W+B_W+\Mm_W)+\bar E$$
for some $\bar E\geq 0$.  By Lemma \ref{lem: special proper log smooth model}, possibly replacing $(\bar X,\bar B,\Mm)$ with a higher model, we may assume that there is a decomposition $\bar B=\bar B^h+\bar B^v$ such that
\begin{itemize}
    \item[(ii)] $\bar B^h\geq 0$ and $\bar B^v$ is reduced,
    \item[(iii)] $\bar B^v$ is vertical over $V$, and
    \item[(iv)] for any $t\in (0,1]$, all lc centers of $(\bar X,\bar B-t\bar B^v,\Mm)$ dominate $V$.
\end{itemize}
Let $\bar D:=g^*D_W+\bar E$. Then $(\bar X,\bar B+\bar D)$ is log smooth and $\bar D\sim_{\Rr,U}K_{\bar X}+\bar B+\Mm_{\bar X}$. Since $\Supp D_W$ contains any vertical lc center of $(W,B_W,\Mm)$, by \cite[Lemma 3.7]{HL21a} we have $\Supp\bar B^v\subset\Supp\bar D$. Thus we may find $\bar G,\bar H\ge0$ and write $\bar D=\bar G+\bar H$ such that
\begin{itemize}
    \item[(v)] $K_{\bar X}+\bar B+\Mm_{\bar X}\sim_{\Rr,U}\bar G+\bar H$,
    \item[(vi)] $\Supp\bar B^v\subset\Supp\bar G\subset\Supp\lfloor\bar B\rfloor$, and
    \item[(vii)] no component of $\bar H$ is contained in $\lfloor\bar B\rfloor$ and $(\bar X,\bar B+\bar H)$ is log smooth.
\end{itemize}
We fix a real number $t_1\in(0,1)$ such that $\bar B-t_1\bar G\geq 0$. For any $t\in (0,t_1]$, by (ii)(iii)(iv)(vi), any lc center of $(\bar X,\bar B-t\bar G,\Mm)$ dominates $V$. By (i)(v) and Lemma \ref{lem: property of numerical and Iitaka dimension}(2), we have $$\kappa_{\sigma}(\bar X/U,K_{\bar X}+\bar B-t\bar G+\Mm_{\bar X})=\dim V-\dim U$$ and $$\kappa_{\sigma}(\bar X/V,K_{\bar X}+\bar B-t\bar G+\Mm_{\bar X})=\kappa_{\iota}(\bar X/V,K_{\bar X}+\bar B-t\bar G+\Mm_{\bar X})=0.$$ 
Then by Proposition \ref{prop: prop 3.4 has19 abu ver} we obtain
\begin{itemize}
    \item[(viii)] $(\bar X,\bar B-t\bar G,\Mm)/U$ has a good minimal model for any $t\in (0,t_1]$.
\end{itemize}

Since $(\bar X,\bar B,\Mm)$ is a log smooth model of $(X,B,\Mm)$, we may run a $(K_{\bar X}+\bar B+\Mm_{\bar X})$-MMP$/X$ which terminates with a dlt model $(Y,B_Y,\Mm)$ of $(X,B,\Mm)$ with induced morphism $f: Y\rightarrow X$ and birational map $\phi: \bar X\dashrightarrow Y$. Let $G_Y$ and $H_Y$ be the strict transforms of $\bar G$ and $\bar H$ on $Y$ respectively. Then $K_Y+B_Y+\Mm_Y\sim_{\Rr,U}G_Y+H_Y$. By (vii) and Lemma \ref{lem: still an mmp under perturbation}, there exists $0<t_2<t_1$ such that $(\bar Y,\bar B+t_2\bar H,\Mm)$ is dlt and $\phi$ is a $(K_{\bar X}+\bar B+t\bar H+\Mm_{\bar X})$-MMP$/X$ as well as a $(K_{\bar X}+\bar B-t\bar G+\Mm_{\bar X})$-MMP$/X$ for any $t\in (0,t_2]$. Then $(Y,B_Y+t_2H_Y,\Mm)$ is dlt, and by (viii) and \cite[Theorem 2.24, Lemma 3.9(2)]{HL21a}, $(Y,B_Y-tG,\Mm))/U$ has a good minimal model for any $t\in (0,t_2]$. $N_\sigma(Y/U,K_Y+B_Y+tH_Y+\Mm_Y)$ is well-defined as a divisor since $K_Y+B_Y+tH_Y+\Mm_Y\sim_{\Rr,U}G_Y+(1+t)H_Y$ is effective for any $t\ge0$.
By Lemma \ref{lem: Has20 2.4 rel ver}, we may pick $0<t_0<t_2$ such that $\Supp N_{\sigma}(Y/U,K_Y+B_Y+tH_Y+\Mm_Y)$ does not depend on $t$ for any $t\in (0,t_0]$.

\begin{center}$\xymatrix{
Y\ar@{->}[dr] & \bar X\ar@{-->}[l]_{\phi}\ar@{->}[r]^{g} & W\ar@{->}[ld]_{h}\ar@{->}[d]^{\psi}\\
& X\ar@{-->}[r] & V 
}$
\end{center}

We may replace $(X,B,\Mm)$ with $(Y,B_Y,\Mm)$ and let $G:=G_Y$ and $H:=H_Y$, and assume that $(X,B,\Mm),G,H$ and $t_0$ satisfy (I)(II)(III). In what follows, we forget all other auxiliary varieties and divisors constructed in this step.

\medskip

\noindent\textbf{Step 2}. For any $t\in (0,t_0]$, by (III.2), $(X,B-\frac{t}{1+t}G,\Mm)/U$ has a good minimal model. Since
$$K_X+B+tH+\Mm_X\sim_{\Rr,U}(1+t)(K_X+B-\frac{t}{1+t}G+\Mm_X),$$
by (III.1) and \cite[Theorem 2.24, Lemma 3.9(2), 4.2]{HL21a}, we may run a $(K_X+B+tH+\Mm_X)$-MMP$/U$ $\phi_t: X\dashrightarrow X_t$ which terminates with a good minimal model $(X_t,B_t+tH_t,\Mm)/U$ of $(X,B+tH,\Mm)/U$. By Lemma \ref{lem: nz for lc divisor}(2), the divisors contracted by $\phi_t$ are precisely the components of $\Supp N_{\sigma}(X/U,K_X+B+tH+\Mm_X)$.
Since $\Supp N_{\sigma}(X/U,K_X+B+tH+\Mm_X)$ does not depend on $t\in(0,t_0]$ by (III.1), each MMP $\phi_t$ contracts precisely the components of $\Supp N_{\sigma}(X/U,K_X+B+t_0H+\Mm_X)$. We let $X_0:=X_{t_0},B_0:=B_{t_0}$, and $H_0:=H_{t_0}$. Then $X_0$ and $X_t$ are isomorphic in codimension $1$, and $K_{X_0}+B_0+\Mm_{X_0}$ is a movable$/U$ $\Rr$-divisor. By the negativity lemma, $(X_t,B_t+tH_t,\Mm)/U$ is a good minimal model of $(X_0,B_0+tH_0,\Mm)/U$ for any $t\in (0,t_0]$. 

\begin{claim}\label{claim: Has20 3.3 rel gpair}
We may run a $(K_{X_0}+B_0+\Mm_{X_0})$-MMP$/U$ with scaling of $H_0$ 
$$(X_0,B_0,\Mm)\dashrightarrow (X_1,B_1,\Mm)\dashrightarrow\dots \dashrightarrow(X_i,B_i,\Mm)\dashrightarrow\dots$$
with scaling numbers 
$$\lambda_i:=\inf\{t\geq 0 \mid K_{X_i}+B_i+tH_i+\Mm_{X_i}\text{ is nef}/U\},$$
where $H_i$ is the strict transform of $H$ on $X_i$, which consists only of flips such that
\begin{enumerate}
    \item either the MMP$/U$ terminates with a minimal model, or $\lim_{i\rightarrow+\infty}\lambda_i=0$,
    \item for any $i\geq 1$ and $\lambda\in [\lambda_i,\lambda_{i-1}]$, $(X_i,B_i+\lambda H_i,\Mm)/U$ is a good minimal model of both $(X,B+\lambda H,\Mm)$ and $(X_0,B_0+\lambda H_0,\Mm)/U$, and
    \item the MMP only contracts sub-varieties of $\Supp\lfloor B_0\rfloor$.
\end{enumerate}
\end{claim}
\begin{proof}
Since $K_{X_0}+B_0+\Mm_{X_0}$ is a movable$/U$ $\Rr$-divisor, by Lemma \ref{lem: limit of movable divisors mmp only contain flips}, any $(K_{X_0}+B_0+\Mm_{X_0})$-MMP$/U$ only contains flips. (1) follows from Lemma \ref{lem: gmmp scaling numbers go to 0}. For any $i\geq 1$ and $\lambda\in [\lambda_i,\lambda_{i-1}]$, $(X_i,B_i+\lambda H_i,\Mm)$ is dlt and $K_{X_i}+B_i+\lambda H_i+\Mm_{X_i}$ is nef$/U$. Since the induced birational maps $X_0\dashrightarrow X_\lambda$ and $X_i\rightarrow X_\lambda$ are both small, by Lemma \ref{lem: hmx18 2.7.3 gpair rel} and \cite[Lemma 3.9(2)]{HL21a}, we get (2).

Let $X_i\rightarrow Z_i\leftarrow X_{i+1}$ be the $i$-th step of the MMP where $X_i\rightarrow Z_i$ the flipping contraction. Then for any flipping curve $C_i$ of $X_i\rightarrow Z_i$, we have $(K_{X_i}+B_i+\Mm_{X_i})\cdot C_i<0$ and $H_i\cdot C_i>0$. Let $G_i$ be the strict transform of $G$ on $X_i$. Then $0>(K_{X_i}+B_i+\Mm_{X_i}-H_i)\cdot C_i=G_i\cdot C_i$. Thus $C_i\subset\Supp G_i$. Since $\Supp G\subset\Supp\lfloor B\rfloor$, $\Supp G_i\subset\Supp\lfloor B_i\rfloor$, and we get (3).
\end{proof}

\begin{claim}\label{claim: Has22a 3.14 step 2 reduction}
Let
$$(X_0,B_0,\Mm)\dashrightarrow (X_1,B_1,\Mm)\dashrightarrow\dots \dashrightarrow (X_i,B_i,\Mm)\dashrightarrow\dots,$$
$\lambda_i$, and $H_i$ be the MMP$/U$, the scaling numbers, and the strict transform of $H$ on $X_i$ for each $i$ as in Claim \ref{claim: Has20 3.3 rel gpair} respectively. If the MMP$/U$ terminates, then Theorem \ref{thm: Has22a 3.14 rel ver} holds.
\end{claim}
\begin{proof}
Let $\lambda_{-1}:=t_0$. If the MMP$/U$ terminates, then $\lambda_{l-1}>\lambda_l=0$ for some $l\in\mathbb N$. By Claim \ref{claim: Has20 3.3 rel gpair}(2), for any $t\in (0,\lambda_{l-1}]$, $K_{X_l}+B_l+tH_l+\Mm_{X_l}$ is nef$/U$, and $a(P,X,B+tH,\Mm)\leq a(P,X_l,B_l+tH_l,\Mm)$ for any prime divisor $P$ on $X$ that is exceptional$/X_l$. Letting $t\rightarrow 0$, we have that  $K_{X_l}+B_l+\Mm_{X_l}$ is nef$/U$ and $a(P,X,B,\Mm)\leq a(P,X_l,B_l,\Mm)$ for any prime divisor $P$ on $X$ that is exceptional$/X_l$. Thus $(X_l,B_l,\Mm)/U$ is a weak lc model of $(X,B,\Mm)/U$. The Claim follows from \cite[Lemma 3.15]{HL21a}.
\end{proof}

\noindent\textit{Proof of Theorem \ref{thm: Has22a 3.14 rel ver} continued}. In the following, we let
$$(X_0,B_0,\Mm)\dashrightarrow (X_1,B_1,\Mm)\dashrightarrow\dots\dashrightarrow(X_i,B_i,\Mm)\dashrightarrow\dots,$$
$\lambda_i$, and $H_i$ be the MMP$/U$, the scaling numbers, and the strict transform of $H$ on $X_i$ for each $i$ as in Claim \ref{claim: Has20 3.3 rel gpair} respectively.

\medskip

\noindent\textbf{Step 3}. For every $i$ and lc center $S_{i}$ of $(X_{i},B_{i},\Mm)$, we let $(S_i,B_{S_i},\Mm^{S_i})$ be the g-pair induced by adjunction
$$K_{S_i}+B_{S_i}+\Mm^{S_i}_{S_i}:=(K_{X_i}+B_i+\Mm_{X_i})|_{S_i},$$
and let $H_{S_i}:=H_i|_{S_i}$. For every lc center $S$ of $(X,B,\Mm)$ we let $(S,B_S,\Mm^S)$ be the g-pair induced by adjunction
$$K_S+B_S+\Mm^S_S:=(K_X+B+\Mm_X)|_S$$
and let $H_S:=H|_S$. Since $X_0\dashrightarrow X_i$ is a $(K_{X_0}+B_0+\Mm_{X_0})$-MMP$/U$, $X_0\dashrightarrow X_i$ is an isomorphism near the generic point of $S_i$, the strict transform $S_0$ of $S_i$ on $X_0$ is an lc center of $(X_0,B_0,\Mm)$, hence also an lc center of $(X_0,B_0+t_0H_0,\Mm)$. By the same argument, $\phi_{t_0}:X\dashrightarrow X_0$ is an isomorphism near the generic point of $S_0$ and the strict transform $S$ of $S_0$ on $X$ is an lc center of $(X,B+t_0H,\Mm)$. Since $(X,B,\Mm)$ and $(X,B+tH,\Mm)$ are both dlt and have the same lc centers, $S$ is also an lc center of $(X,B,\Mm)$. In particular, the induced maps $\phi^S_i:S\dashrightarrow S_i$ and $\phi^S_{j,i}:S_j\dashrightarrow S_i$ are birational for any $j\le i$. 

By Lemma \ref{lem: special termination reduce to flip lemma}(1), we may find $m\gg0$ such that $X_m\dashrightarrow X_i$ is an isomorphism near the generic point of any lc center $S_m$ of $(X_m,B_m,\Mm)$ and any $i\ge m$. By Lemma \ref{lem: special termination reduce to flip lemma}(2.a), possibly replacing $m$, we may assume that the induced $\phi^S_{m,i}:S_m\to S_i$ is small for any lc center $S_m$ of $(X_m,B_m,\Mm)$ and any $i\geq m$.

Then we only need to show that for any lc center $S_m$ of $(X_m,B_m,\Mm)$ of dimension $\geq 1$, $(S_m,B_{S_m},\Mm^S)/U$ has a log minimal model. Indeed, if this is the case, then by Lemma \ref{lem: special termination reduce to flip lemma}, \cite[Remark 3.25, Theorem 4.1]{HL22}, Claim \ref{claim: Has20 3.3 rel gpair}(1), and induction on the dimension of lc centers, the $(K_{X_0}+B_0+\Mm_{X_0})$-MMP$/U$ above will induce isomorphisms on any lc center $S_m$ of $(X_m,B_m,\Mm)$ for $m\gg0$. But by Claim \ref{claim: Has20 3.3 rel gpair}(3) the $(K_{X_0}+B_0+\Mm_{X_0})$-MMP$/U$ only contracts sub-varieties of $\Supp\lfloor B_0\rfloor$, so it must terminate.



\medskip\noindent\textbf{Step 4}. We prove the following claim.

\begin{claim}\label{claim: Has22a 3.14 step 4 abcd}
There exists a $\Qq$-factorial lc g-pair $(T,B_T,\Mm^S)/U$ and a birational morphism $\psi: T\rightarrow S_m$ satisfying the following:
\begin{enumerate}
    \item For any prime divisor $D$ on $S$ such that $a(D,S_m,B_{S_m},\Mm^S)<a(D,S,B_S,\Mm^S)$, $D$ is on $T$ and is a $\psi$-exceptional.
    \item $$B_T=\sum_{D\text{ is a prime divisor on }T}(1-a(D,S,B_S,\Mm^S))D.$$
    \item For any $i\geq m$ and any prime divisor $Q$ over $S$, we have $$a(Q,S,B_S+\lambda_iH_S,\Mm^S)\leq a(Q,S_i,B_{S_i}+\lambda_iH_{S_i},\Mm^S).$$
    \item For any prime divisor $Q'$ over $S_m$, we have
    $$a(Q',S_m,B_{S_m},\Mm^S)\leq a(Q',T,B_T,\Mm^S).$$
\end{enumerate}
\end{claim}
\begin{proof}
By Claim \ref{claim: Has20 3.3 rel gpair}(2), $(X_i,B_i+\lambda_iH_i,\Mm)/U$ is a good minimal model of $(X,B+\lambda_iH,\Mm)$, hence (3) holds.

Since $\phi_i$ does not extract any divisor, $a(P,X_i,B_i,\Mm)\leq a(P,X,B,\Mm)$ for any prime divisor $P$ on $X_i$. Since  $\sigma_P(X/U,K_X+B+\Mm_X)=0$ for any prime divisor $P$ over $X$ such that $a(P,X,B,\Mm)<1$ and $\Center_XP\cap\Nklt(X,B,\Mm)\not=\emptyset$, by Lemma \ref{lem: Has22a 3.5 rel} and since $\phi_{m,i}$ is small for any $i\geq m$, $a(D,S_i,B_{S_i},\Mm^S)\leq a(D,S,B_S,\Mm^S)$ for any prime divisor $D$ on $S_m$ and $i\geq m$. Thus  $a(D,S_i,B_{S_i}+\lambda_iH_{S_i},\Mm^S)\leq a(D,S,B_S,\Mm^S)$ for any prime divisor $D$ on $S_m$ and $i\geq m$. By (3), for any $i\geq m$ we have
\begin{align*}
    a(D,S,B_S+\lambda_iH_S,\Mm^S)&\leq a(D,S_i,B_{S_i}+\lambda_iH_{S_i},\Mm^S)\\
    &=a(D,S_m,B_{S_m}+\lambda_iH_{S_m},\Mm^S)\leq a(D,S,B_S,\Mm^S).
\end{align*}

Letting $i\rightarrow+\infty$, we have
$$a(D,S,B_S,\Mm^S)=a(D,S_m,B_{S_m},\Mm^S)$$
for any prime divisor $D$ on $S_m$. We define
$$\mathcal{C}:=\{D\mid D\text{ is a prime divisor on }S, a(D,S_m,B_{S_m},\Mm^S)<a(D,S,B_S,\Mm^S)\}.$$
Then any element of $\mathcal{C}$ is exceptional over $S_m$. Thus for any $D\in\mathcal{C}$, by (3), we have
\begin{align*}
    a(D,S,B_S+\lambda_mH_S,\Mm^S)&\leq a(D,S_m,B_{S_m}+\lambda_mH_{S_m},\Mm^S)\\
    &\leq a(D,S_m,B_{S_m},\Mm^S)<a(D,S,B_S,\Mm^S)\leq 1.
\end{align*}
Since any element of $\mathcal{C}$ is a prime divisor on $S$, any element of $\mathcal{C}$ is a component of $H_S$. Thus $\mathcal{C}$ is a finite set, and for every $D\in\mathcal{C}$, since $\lambda_m<t_0$, we have
\begin{align*}
  0\leq &a(D,S,B_S+t_0H_S,\Mm^S)<a(D,S,B_S+\lambda_mH_S,\Mm^S)\\ \leq &a(D,S_m,B_{S_m}+\lambda_mH_{S_m},\Mm^S)\leq  a(D,S_m,B_{S_m},\Mm^S)<a(D,S,B_S,\Mm^S)\leq 1.  
\end{align*}
Thus $0<a(D,S_m,B_{S_m},\Mm^S)<1$ for any $D\in\mathcal{C}$. By \cite[Lemma 3.4]{Has22a}, there exists a birational morphism $\psi: T\rightarrow S_m$ from a $\Qq$-factorial variety $T$ which extracts exactly divisors contained in $\mathcal{C}$. (1) follows immediately from the construction of $\mathcal{C}$. Since $(S,B_S,\Mm^S)$ is lc, there are only finitely many divisors $D$ on $T$ such that $a(D,S,B_S,\Mm^S)\not=1$, hence $B_T\geq 0$ is well-defined, and we get (2).

For any prime divisor $D$ on $T$, if $D$ is $\psi$-exceptional, then $$a(D,S_m,B_{S_m},\Mm^S)< a(D,S,B_S,\Mm^S)\leq 1$$ as $D\in\mathcal{C}$, and if $D$ is not $\psi$-exceptional, then $\Center_{S_m}D$ is a divisor, hence $a(D,S,B_S,\Mm^S)=a(D,S_m,B_{S_m},\Mm^S)\leq 1$. In either case,
$$a(D,S_m,B_{S_m},\Mm^S)\leq a(D,S,B_S,\Mm^S)\leq 1.$$
Since $T$ is $\Qq$-factorial, $K_T+B_T+\Mm^S_T$ is $\Rr$-Cartier, and
$$K_T+B_T+\Mm^S_T\leq\psi^*(K_{S_m}+B_{S_m}+\Mm^S_{S_m}).$$
Thus
$$0\leq a(Q',S_m,B_{S_m},\Mm^S)\leq a(Q',T,B_T,\Mm^S)$$
for any prime divisor $Q'$ over $S_m$, and we get (4). In particular, $(T,B_T,\Mm^S)$ is lc, and the proof is concluded.
\end{proof}

\noindent\textbf{Step 5}. In this step we show that $(T,B_T,\Mm^S)/U$ has a log minimal model. We first prove the following claim:

\begin{claim}\label{lem: Has22a 3.14 claim rel version}
For any prime divisor $D$ over $S$,
\begin{enumerate}
    \item if $D$ is on $S$, then $a(D,S,B_S,\Mm^S)\leq a(D,T,B_T,\Mm^S)$, and
    \item if $D$ is on $T$, then $a(D,T,B_T,\Mm^S)= a(D,S,B_S,\Mm^S)$.
\end{enumerate}
\end{claim}
\begin{proof}
By Claim \ref{claim: Has22a 3.14 step 4 abcd}(2), we only need to show that for any prime divisor $D$ on $S$ that is exceptional over $T$, $a(D,S,B_S,\Mm^S)\leq a(D,T,B_T,\Mm^S)$. By Claim \ref{claim: Has22a 3.14 step 4 abcd}(1)(4), $$a(D,S,B_S,\Mm^S)\leq a(D,S_m,B_m,\Mm^S)\leq a(D,T,B_T,\Mm^S),$$
and we get (1).
\end{proof}

By our assumption, $(S,B_S,\Mm^S)/U$ is a log minimal model of itself, then by Lemma \ref{lem: Has22a 3.9 rel} $(T,B_T,\Mm^S)/U$ also has a log minimal model.

\medskip

\noindent\textbf{Step 6}. We conclude the proof in this step. Recall that we only need to show that $(S_m,B_{S_m},\Mm^S)/U$ has a log minimal model.

For any $i\geq m$, since $K_{X_i}+B_i+\lambda_iH_i+\Mm_{X_i}$ is nef$/U$, $K_{S_i}+B_{S_i}+\lambda_iH_{S_i}+\Mm^S_{S_i}=(K_{X_i}+B_i+\lambda_iH_i+\Mm_{X_i})|_{S_i}$ is nef$/U$. Since $\phi^S_{m,i}$ is small, 
$(S_i,B_{S_i}+\lambda_iH_{S_i},\Mm^S)/U$ is a weak lc model of $(S_m,B_{S_m}+\lambda_iH_{S_m},\Mm^S)/U$. Let $h_m: W\rightarrow S_m$ and $h_i: W\rightarrow S_i$ be a resolution of indeterminacies of $\phi^S_{m,i}$. By Lemmas \ref{lem: nz pre properties}(1), \ref{lem: nz keep under pullback}(2), \ref{lem: nz basic properties}(3) and the negativity lemma, for any prime divisor $D$ on $T$ we have 
\begin{align*}
    0&\leq a(D,S_i,B_{S_i}+\lambda_iH_{S_i},\Mm^S)-a(D,S_m,B_{S_m}+\lambda_iH_{S_m},\Mm^S)\\
    &=\sigma_D(S_m/U,K_{S_m}+B_{S_m}+\lambda_iH_{S_m}+\Mm^S_{S_m}).
\end{align*}
By Claim \ref{claim: Has22a 3.14 step 4 abcd}(3), we have
$$\sigma_D(S_m/U,K_{S_m}+B_{S_m}+\lambda_iH_{S_m}+\Mm^S_{S_m})\geq  a(D,S,B_{S}+\lambda_iH_{S},\Mm^S)-a(D,S_m,B_{S_m}+\lambda_iH_{S_m},\Mm^S).$$
By Claim \ref{claim: Has22a 3.14 step 4 abcd}(2), $a(D,S,B_S,\Mm^S)=a(D,T,B_T,\Mm^S)$. By Lemma \ref{lem: nz basic properties}(2) and Claims \ref{claim: Has20 3.3 rel gpair}(1) and \ref{claim: Has22a 3.14 step 4 abcd}(4), for any prime divisor $D$ on $T$,
\begin{align*}
    &\sigma_D(S_m/U,K_{S_m}+B_{S_m}+\Mm^S_{S_m})\\
    =&\lim_{i\rightarrow+\infty}\sigma_D(S_m/U,K_{S_m}+B_{S_m}+\lambda_iH_{S_m}+\Mm^S_{S_m})\\
    \geq&\lim_{i\rightarrow+\infty}(a(D,S,B_S+\lambda_iH_S,\Mm^S)-a(D,S_m,B_{S_m}+\lambda_iH_{S_m},\Mm^S))\\
    =&a(D,S,B_S,\Mm^S)-a(D,S_m,B_{S_m},\Mm^S)\\
    =&a(D,T,B_T,\Mm^T)-a(D,S_m,B_{S_m},\Mm^S)\geq 0.
\end{align*}
Since $(T,B_T,\Mm^S)/U$ has a log minimal model by Step 5, by Lemma \ref{lem: Has22a 3.6 rel}, $(S_m,B_{S_m},\Mm^S)/U$ has a log minimal model, and we are done.
\end{proof}

\begin{thm}[{cf. \cite[Theorem 3.15]{Has22a}}]\label{thm: Has22a 3.15 rel ver}
Let $(X,B,\Mm)/U$ be a $\Qq$-factorial NQC dlt g-pair and $A\geq 0$ an $\mathbb{R}$-divisor on $X$ such that $(X,B+A,\Mm)/U$ is lc and $K_{X}+B+\Mm_{X}+A$ is nef$/U$. Then for any $(K_{X}+B+\Mm_{X})$-MMP$/U$ with scaling of $A$
$$(X,B, \Mm)=:(X_{0},B_{0},\Mm) \dashrightarrow (X_{1}, B_{1},\Mm) \dashrightarrow \dots \dashrightarrow (X_{i}, B_{i},\Mm) \dashrightarrow \dots,$$
with scaling numbers
$$\lambda_i:=\inf\{t\geq 0 \mid K_{X_i}+B_i+tA_i+\Mm_{X_i}\text{ is nef}/U\},$$
where $A_i$ is the strict transform of $A$ on $X_i$, if $\lambda_i>0$ for each $i$ and $\lim_{i\rightarrow+\infty}\lambda_i=0$, then there are only finitely many $i$ such that $K_{X_i}+B_i+\Mm_{X_i}$ is log abundant$/U$ with respect to $(X_i,B_i,\Mm)$.
\end{thm}
\begin{proof} We apply induction on the dimension. Suppose that the theorem holds in dimension $\leq \dim X-1$ but the theorem does not hold. Then there exists a $(K_{X}+B+\Mm_{X})$-MMP$/U$ with scaling of $A$ as in the statement of the theorem such that  $K_{X_i}+B_i+\Mm_{X_i}$ is log abundant$/U$ with respect to $(X_i,B_i,\Mm)$ for infinitely many $i$. Let $\phi_{i,j}: X_i\dashrightarrow X_j$ be the induced birational map. Possibly replacing $(X,B,\Mm)$ with $(X_m,B_m,\Mm)$ for some $m\gg 0$, we may assume that the maps $\phi_{i,j}$ are small for any $i,j$. 

We prove the following claim.

\begin{claim}\label{claim: has20 3.4 step 2 gpair rel}
If there exists $m\gg 0$ such that $\phi_{m,i}|_S$ is an isomorphism for any lc center $S$ of $(X_m,B_m,\Mm)$ and $i\geq m$, then Theorem \ref{thm: Has22a 3.15 rel ver} holds.
\end{claim}
\begin{proof}
Possibly replacing $(X,B,\Mm)$ with $(X_m,B_m,\Mm)$ we may assume that $\phi_{i,j}|_{\Nklt(X_i,B_i,\Mm)}$ is an isomorphism for any $i,j$ and $K_X+B+\Mm_X$ is abundant$/U$. Since $\lim_{i\rightarrow+\infty}\lambda_i=0$ and $\phi_{i,j}$ are small for any $i,j$, $K_X+B+\Mm_X$ is a movable$/U$ $\Rr$-divisor, hence $K_X+B+\Mm_X$ is pseudo-effective$/U$. Notice that $(X_i,B_i,\Mm)$ are all dlt and $\Qq$-factorial. Let $D$ be a component of $\lf B_i\rf$. Then $\phi_{i,i+1}|_D$ being an isomorphism implies that the flip $\phi_{i,i+1}$ is an isomorphism near $D$. Therefore $\phi_{i,i+1}$ is an isomorphism on a neighborhood of $\lf B_i\rf$. By Lemma \ref{lem: Has22a 3.8 rel}, ${\bf B}_{-}(K_{X}+B+\Mm_{X}/U)$ does not intersect $\Supp\lfloor B\rfloor$, and $\sigma_P(X/U,K_X+B+\Mm_X)=0$ for any prime divisor $P$ over $X$ such that $\Center_XP\cap\Supp\lfloor B\rfloor\not=\emptyset$. In particular, $(K_X+B+\Mm_X)|_S$ is nef$/U$ for any lc center $S$ of $(X,B,\Mm)$. By Theorem \ref{thm: Has22a 3.14 rel ver}, $(X,B,\Mm)/U$ has a log minimal model, but this contradicts \cite[Theorem 4.1]{HL22} so we are done.
\end{proof}

\noindent\textit{Proof of Theorem \ref{thm: Has22a 3.15 rel ver} continued}. We let $\phi_i:=\phi_{i,i+1}$ for every $i$ and $X_i\rightarrow Z_i\leftarrow X_{i+1}$ the flip defined by $\phi_i$. By Claim \ref{claim: has20 3.4 step 2 gpair rel}, we only need to show that for any lc center $S$ of $(X,B,\Mm)$, the MMP terminates along $S$ after finitely many steps. By induction on the dimension of lc centers, we may assume that $\phi_i$ induces an isomorphism for every $k$-dimensional lc centers and $i\gg 0$, where $k<d=\dim S$. Let $S_i$ be the strict transform of $S$ on $X_i$ and $(S_i,B_{S_i},\Mm^{S})$ the g-pair given by adjunction
$$K_{S_i}+B_{S_i}+\Mm_{S}^{S_i}:=(K_{X_i}+B_i+\Mm_{X_i})|_{S_i}.$$
Let $(S_i',B_{S_i'},\Mm^S)$ be a dlt model of $(S_i,B_{S_i},\Mm^{S})$. By Lemma \ref{lem: special termination reduce to flip lemma}, for $i\gg0$, the $(K_X+B+\Mm_X)$-MMP$/U$ with scaling of $A$ induces a $(K_{S'_i}+B_{S'_i}+\Mm_{S}^{S'_i})$-MMP$/T$ with scaling of $A_{S_i'}$ such that the limit of the scaling numbers is $0$, where $A_{S_i'}$ is the pullback of $A_i$ on $S_i'$ and $T$ is the normalization of the image of $S_i$ in $U$. Since $K_{X_j}+B_j+\Mm_{X_j}$ is log abundant$/U$ with respect to $(X_j,B_j,\Mm)$ for infinitely many $j$, $K_{S_j'}+B_{S_j'}+\Mm^S_{S'_j}$ is log abundant$/T$ with respect to $(S_j',B_{S_j'},\Mm^S)$ for infinitely many $j$. By Theorem \ref{thm: Has22a 3.15 rel ver} in dimension $<\dim X$, the $(K_{S'_j}+B_{S'_j}+\Mm_{S'_j}^{S})$-MMP$/T$ terminates, i.e the top horizontal maps in Lemma \ref{lem: special termination reduce to flip lemma}(3) are isomorphisms for $i\gg0$. Therefore the bottom horizontal maps in Lemma \ref{lem: special termination reduce to flip lemma}(3) must also be isomorphisms for $i\gg0$ by the contructions in the proof. Thus there exists $m\gg 0$ such that $\phi_{m,i}|_S$ is an isomorphism for any lc center $S_m$ of $(X_m,B_m,\Mm)$ and $i\geq m$ and we are done.
\end{proof}

\section{Log abundance under the MMP}

This section is similar to \cite[Section 3 and Theorem 4.1]{Has22b}.

\begin{thm}[{cf. \cite[Theorem 3.5]{Has22b}}]\label{thm: Has22b 3.5 rel ver}
Let $(X,B,\Mm)/U$ be an NQC lc g-pair and $\pi: X\rightarrow Z$ a projective morphism$/U$ such that $Z$ is normal quasi-projective. Let $C\geq 0$ be an $\Rr$-divisor on $X$, $A_Z$ an ample$/U$ $\Rr$-divisor on $Z$, and $0\leq A\sim_{\Rr,U}\pi^*A_Z$ an $\Rr$-divisor on $X$, such that
\begin{enumerate}
    \item $C$ does not contain any lc center of $(X,B,\Mm)$,
    \item $K_X+B+C+\Mm_X\sim_{\Rr,Z}0$, and
    \item $(X,B+A,\Mm)$ is lc.
\end{enumerate}
Then $K_X+B+A+\Mm_X$ is abundant$/U$.
\end{thm}
\begin{proof}
Possibly replacing $\pi$ with the contraction in the Stein factorization of $\pi$, we may assume that $\pi$ is a contraction. Possibly replacing $Z\rightarrow U$ with the Stein factorization of $Z\rightarrow U$, we may assume that $Z\rightarrow U$ is a contraction. Let $F$ be a very general fiber of $X\rightarrow U$ and $F_Z:=\pi(F)$. Then $F_Z$ is a very general fiber of $Z\rightarrow U$. Possibly replacing $(X,B,\Mm),A,C,Z,A_Z,\pi,U$ with $(F,B|_F,\Mm|_F),A|_F,C|_F,F_Z,A_Z|_{F_Z},\pi|_F,\{pt\}$, we may assume that $U=\{pt\}$. The theorem follows from \cite[Theorem 3.5]{Has22b}.  Note that we remove the $\Rr$-Cartier assumption of $C$ as it is immediate from (2).
 \end{proof}
 
\begin{lem}[{cf. \cite[Lemma 3.6]{Has22b}}]\label{lem: Has22b 3.6 rel ver}
Let $(X,B,\Mm)/U$ be an NQC lc g-pair and $\pi: X\rightarrow Z$ a projective morphism$/U$ such that $Z$ is normal quasi-projective. Let $C\geq 0$ be an $\Rr$-divisor on $X$, $A_Z$ an ample$/U$ $\Rr$-divisor on $Z$, and $0\leq A\sim_{\Rr,U}\pi^*A_Z$ an $\Rr$-divisor on $X$, such that
\begin{enumerate}
    \item $C$ does not contain any lc center of $(X,B,\Mm)$,
    \item $K_X+B+C+\Mm_X\sim_{\Rr,Z}0$, and
    \item $(X,B+A,\Mm)$ is lc.
\end{enumerate}
Let $h: W\rightarrow X$ be a log resolution of $(X,\Supp B)$ such that $\Mm$ descends to $W$, and $B_W\geq 0$ an $\Rr$-divisor on $W$ such that $(W,B_W+h^*A)$ is lc and $(K_W+B_W+\Mm_W-h^*(K_X+B+\Mm_X))^{\geq 0}$ is $h$-exceptional. Then $K_W+B_W+h^*A+\Mm_W$ is abundant$/U$. 
\end{lem}
\begin{proof}
Possibly replacing $\pi$ with the contraction in the Stein factorization of $\pi$, we may assume that $\pi$ is a contraction. Possibly replacing $Z\rightarrow U$ with the Stein factorization of $Z\rightarrow U$, we may assume that $Z\rightarrow U$ is a contraction. Let $F_w$ be a very general fiber of $W\rightarrow U$, $F:=h(F_W)$, and $F_Z:=\pi(F)$. Then $F$ and $F_Z$ are very general fibers of $X\rightarrow U$ and $Z\rightarrow U$ respectively. Possibly replacing $(X,B,\Mm),A,C,Z,A_Z,\pi,U,W,h,B_W$ with $(F,B|_F,\Mm|_F),A|_F,C|_F,F_Z,A_Z|_{F_Z},\pi|_F,\{pt\},F_W,h|_{F_W},B_W|_{F_W}$, we may assume that $U=\{pt\}$. The theorem follows from \cite[Theorem 3.5, Lemma 3.6]{Has22b}. Note that we remove the $\Rr$-Cartier assumption of $C$ as it is immediate from (2).
\end{proof}

\begin{thm}[{cf. \cite[Theorem 4.1]{Has22b}}]\label{thm: Has22b 4.1 rel ver}
Let $(X,B,\Mm)/U$ be an NQC lc g-pair and $\pi: X\rightarrow Z$ a projective morphism$/U$ such that $Z$ is normal quasi-projective. Let $C\geq 0$ be an $\Rr$-divisor on $X$, $A_Z$ an ample$/U$ $\Rr$-divisor on $Z$, and $0\leq A\sim_{\Rr,U}\pi^*A_Z$ an $\Rr$-divisor on $X$, such that
\begin{enumerate}
    \item $C$ does not contain any lc center of $(X,B,\Mm)$,
    \item $K_X+B+C+\Mm_X\sim_{\Rr,Z}0$, and
    \item $(X,\Delta:=B+A,\Mm)$ is lc and $\Nklt(X,B,\Mm)=\Nklt(X,\Delta,\Mm)$.
\end{enumerate}
Then for any $(K_{X}+\Delta+\Mm_{X})$-MMP$/U$
$$(X,\Delta,\Mm):=(X_0,\Delta_0,\Mm)\dashrightarrow (X_1,\Delta_1,\Mm)\dashrightarrow\dots\dashrightarrow (X_i,\Delta_i,\Mm)\dashrightarrow\dots,$$
$K_{X_{i}}+\Delta_{i}+\Mm_{X_{i}}$ is log abundant$/U$ with respect to $(X_{i},\Delta_{i},\Mm)$ for every $i$. 
\end{thm}

\begin{proof}
For each $i$, we let $\phi_i: X\dashrightarrow X_i$ be the induced birational map. 

By Theorem \ref{thm: Has22b 3.5 rel ver}, $K_{X}+B+A+\Mm_{X}$ is abundant$/U$. By Lemma \ref{lem: property of numerical and Iitaka dimension}(6), $K_{X_{i}}+\Delta_{i}+\Mm_{X_{i}}$ is abundant$/U$ for any $i$. Thus we only need to prove that $(K_{X_i}+\Delta_i+\Mm_{X_i})|_{S_i}$ is abundant$/U$ for any lc center $S_i$ of $(X_i,\Delta_i,\Mm)$.

Fix $i$ and an lc center $S_i$ of $(X_i,\Delta_i,\Mm)$. Then there exists an lc center $S$ of $(X,\Delta,\Mm)$ such that $\phi_i|_S: S\dashrightarrow S_i$ is a birational map. Let $(X',B',\Mm)$ be a dlt model of $(X,B,\Mm)$ with induced birational morphism $f: X'\rightarrow X$ such that there exists a component $S'$ of $\lfloor\Delta'\rfloor$ which dominates $S$. Let $C':=f^*C,A':=f^*A$, and $\Delta':=B'+A'$. Since $\Nklt(X,B,\Mm)=f(\Nklt(X',\Delta',\Mm))$, $(X',\Delta',\Mm)$ is a dlt model of $(X,\Delta,\Mm)$. By Lemma \ref{lem: lift mmp}, we may run a $(K_{X'}+\Delta'+\Mm_{X'})$-MMP$/U$ and get a dlt model $(X_i',\Delta_i',\Mm)$ of $(X_i,\Delta_i,\Mm)$ such that the strict transform $S'_i$ of $S'$ on $X_i'$ is a component of $\lfloor\Delta_i'\rfloor$. Then $(K_{X_i}+\Delta_i+\Mm_{X_i})|_{S_i}$ is abundant$/U$ if and only if $(K_{X'_i}+\Delta'_i+\Mm_{X'_i})|_{S'_i}$ is abundant$/U$. Moreover, we have
\begin{itemize}
    \item $C'$ does not contain any lc center of $(X',B',\Mm)$,
    \item $K_{X'}+B'+C'+\Mm_{X'}\sim_{\Rr,Z}0$,
    \item $(X',\Delta',\Mm)$ is lc, and
    \item  $\Nklt(X',B',\Mm)=\Nklt(X',\Delta',\Mm)$. 
\end{itemize} 
Thus possibly replacing $(X,\Delta,\Mm)\dashrightarrow (X_i,\Delta_i,\Mm)$ with $(X',\Delta',\Mm)\dashrightarrow (X'_i,\Delta',\Mm)$ and $A,B,C$ with $A',B',C'$, we may assume that $(X,\Delta,\Mm)$ is $\Qq$-factorial dlt, $S$ is a component of $\lfloor\Delta\rfloor=\lfloor B\rfloor$, and $S_i$ is a component of $\lfloor\Delta_i\rfloor=\lfloor B_i\rfloor$.

Let $(S,B_S,\Mm^S)/U$ and $(S_i,B_{S_i},\Mm^S)/U$ be the dlt g-pairs induced by the adjunction formulas
$$K_S+B_S+\Mm^S:=(K_X+B+\Mm_X)|_S$$
and
$$K_{S_i}+B_{S_i}+\Mm_{S_i}:=(K_{X_i}+B_i+\Mm_{X_i})|_{S_i}.$$
Let $p: W\rightarrow S$ and $q: W\rightarrow S_i$ be a resolution of indeterminacies of the induced birational map $S\dashrightarrow S_i$ such that $\Mm^S$ descends to $W$, $p$ is a log resolution of $(S,\Supp B_S)$, and $q$ is a log resolution of $(S_i,\Supp B_{S_i})$. Since $A$ is semi-ample$/U$, possibly replacing $A$ with a general member of $|A/U|_{\Rr}$\footnote{The ``general member" of the $\Rr$-linear system $|A/U|_{\Rr}$ is constructed in the following way: we write $A=\sum r_iA_i$ where $r_i\in (0,1)$ are real numbers and $A_i$ are base-point-free$/U$ Cartier divisors. We replace $A$ with $\sum r_iH_i$ where $H_i\in |A_i/U|$ are general members.} and setting $A_S:=A|_S$ and $A_{S_i}:=((\phi_i)_*A)|_{S_i}$, we may assume that 
\begin{itemize}
\item $A_S\geq 0, A_{S_i}\geq 0$,
    \item $(S,\Delta_S:=B_S+A_S,\Mm)$ and $(S_i,\Delta_{S_i}:=B_{S_i}+A_{S_i},\Mm)$ are dlt, and
    \item $p$ is a log resolution of $(S,\Supp\Delta_S)$ and $q$ is a log resolution of $(S_i,\Supp\Delta_{S_i})$.
\end{itemize}
Moreover, since $A$ is general in $|A/U|_{\Rr}$, $p^*A_S=p^{-1}_*A_S$, hence $A_W:=p^*A\leq q^*A_{S_i}$. We may write
$$K_W+B_W+A_W+\Mm^S_W=q^*(K_{S_i}+\Delta_{S_i}+\Mm^S_{S_i})+E$$
and let $\Delta_W:=B_W+A_W$, such that $B_W\geq 0$, $E\geq 0$, and $\Delta_W\wedge E=0$. Then $(W,\Delta_W)$ is log smooth. We may write
$$K_W+B_W+\Mm^S_W=p^*(K_S+\Delta_S+\Mm^S_S)+G_+-G_-,$$
where $G_+\geq 0, G_-\geq 0$, and $G_+\wedge G_-=0$. 

For any component $D$ of $G_+$, 
$$a(D,S,\Delta_S,\Mm^S)>a(D,W,\Delta_W,\Mm^S)=\min\{a(D,S_i,\Delta_{S_i},\Mm^S),1\}.$$
Since $\phi_i$ is $(K_X+\Delta+\Mm_X)$-non-positive, by \cite[Lemma 4.2.10]{Fuj07}, $$a(D,S_i,\Delta_{S_i},\Mm^S)\leq a(D,S,\Delta_{S},\Mm^S).$$
Thus $a(D,S,\Delta_S,\Mm^S)>1$, hence $D$ is $p$-exceptional (since any divisor on $S$ has log discrepancy $\le 1$.). We conclude that $G_+$ is $p$-exceptional.

\begin{center}$\xymatrix{
W\ar@{->}[d]_{p}\ar@{->}[dr]^{q}&  \\
S\ar@{-->}[r]_{\phi_i|_S}\ar@{->}[d]_{\pi_S}& S_i \\
Z\ar@{->}[r]& U 
}$
\end{center}

Let $\pi_S:=\pi|_S$ and $C_S:=C|_S$. Since $C$ does not contain any lc center of $(X,B,\Mm)$, $S$ is not a component of $C$, hence $C_S\geq 0$. Then $(S,B_S,\Mm^S)/U$ is an NQC lc g-pair, $\pi_S: S\rightarrow Z$ is a projective morphism$/U$, $Z$ is normal quasi-projective, $C_S\geq 0$ is an $\Rr$ divisor on $X$, $0\leq A_S\sim_{\Rr,U}\pi_S^*A_Z$, such that
\begin{itemize}
    \item $C_S$ does not contain any lc center of $(S,B_S,\Mm^S)$,
    \item $K_S+B_S+C_S+\Mm^S_S\sim_{\Rr,Z}0$,
    \item $(S,\Delta_S=B_S+A_S,\Mm^S)$ is lc and $\Nklt(S,B_S,\Mm^S)=\Nklt(S,\Delta_S,\Mm^S)$,
    \item $p: W\rightarrow S$ is a log resolution of $(S,\Supp B_S)$ such that $\Mm^S$ descends to $W$, $B_W\geq 0$ is an $\Rr$-divisor on $W$ such that $(W,B_W+p^*A_S)$ is lc and $$
    (K_W+B_W+\Mm^S_W-p^*(K_S+B_S+\Mm^S_S))^{\geq 0}=G_+$$ is $p$-exceptional. 
\end{itemize}
By Lemma \ref{lem: Has22b 3.6 rel ver}, $K_W+\Delta_W+\Mm^S_W$ is abundant$/U$. By Lemma \ref{lem: property of numerical and Iitaka dimension}(3), $K_{S_i}+\Delta_{S_i}+\Mm^S_{S_i}=(K_{X_i}+\Delta_i+\Mm_{X_i})|_{S_i}$ is abundant$/U$ and we are done.
\end{proof}

\begin{thm}[{cf. \cite[Theroem 1.3]{Has22b}}]\label{thm: Has22b 1.3 rel ver}
Let $(X,B,\Mm)/U$ be an NQC lc g-pair and $A\geq 0$ an ample$/U$ $\Rr$-divisor such that $(X,\Delta:=B+A,\Mm)$ is lc and $\Nklt(X,B,\Mm)=\Nklt(X,B+A,\Mm)$. Let $(Y,\Delta_Y,\Mm)$ be a dlt model of $(X,\Delta,\Mm)$. Then for any partial $(K_Y+\Delta_Y+\Mm_Y)$-MMP$/U$
$$\phi:(Y,\Delta_Y,\Mm)\dashrightarrow (Y',\Delta_Y',\Mm),$$
$K_{Y'}+\Delta_Y'+\Mm_{Y'}$ is log abundant$/U$ with respect to $(Y',\Delta_Y',\Mm)$.
\end{thm}
\begin{proof}
It follows from Theorem \ref{thm: Has22b 4.1 rel ver}.
\end{proof}

\section{Proof of the main theorems}

\begin{proof}[Proof of Theorem \ref{thm: Has22a 3.17 rel ver intro}]
By \cite[Lemma 4.2]{HL21a}, possibly replacing $A$ by a sufficiently general member, we may assume that $\Nklt(X,B,\Mm)=\Nklt(X,\Delta,\Mm)$.

First we prove (1). Let $(Y,\Delta_Y,\Mm)$ be a dlt model of $(X,\Delta,\Mm)$. By \cite[Theorem 3.14]{HL21a}, we only need to show that $(Y,\Delta_Y,\Mm)/U$ has a log minimal model. We run a $(K_Y+\Delta_Y+\Mm_Y)$-MMP$/U$ 
$$(Y,\Delta_Y,\Mm):=(Y_0,\Delta_0,\Mm)\dashrightarrow(Y_1,\Delta_1,\Mm)\dashrightarrow\dots\dashrightarrow (Y_i,\Delta_i,\Mm)\dashrightarrow\dots$$
with scaling of a general ample$/U$ divisor $H\geq 0$, and let
$$\lambda_i:=\inf\{t\mid t\geq 0, K_{Y_i}+\Delta_i+\lambda_iH_i+\Mm_{Y_i}\text{ is nef}/U\}$$
be the scaling numbers. If $\lambda_i=0$ for some $i$ then $(Y_i,\Delta_i,\Mm)/U$ is a log minimal model of $(Y,\Delta,\Mm)$ and we are done. Thus we may assume that $\lambda_i>0$ for each $i$. By \cite[Theorem 2.24]{HL21a}, $\lim_{i\rightarrow+\infty}\lambda_i=0$. By Theorem \ref{thm: Has22b 1.3 rel ver}, $K_{Y_i}+\Delta_i+\Mm_{Y_i}$ is log abundant$/U$ for each $i$, which contradicts Theorem \ref{thm: Has22a 3.15 rel ver}.

Now we prove (2). We may pick $0<\epsilon\ll 1$ such that $\frac{1}{2}A+\epsilon\Mm_X$ is ample$/U$. Possibly replacing $(X,B,\Mm)$ with $(X,B,(1-\epsilon)\Mm)$ and $A$ with a general member in $|A+\epsilon\Mm_X|_{\Rr}$, we may assume that $\Nklt(X,B)=\Nklt(X,B,\Mm)=\Nklt(X,\Delta,\Mm)$. By \cite[Lemma 5.18]{HL21a}, there exists a birational morphism $h: W\rightarrow X$ such that $\Mm$ descends to $W$ and $\Supp(h^*\Mm_X-\Mm_W)=\Exc(h)$. Let $F:=h^*\Mm_X-\Mm_W$, then $F\geq 0$ and $F$ is exceptional over $X$. In particular, $\Supp F$ does not contain any lc place of $(X,B)$. Thus we may pick $E\geq 0$ on $W$ such that $-E$ is ample$/X$.

Let $K_W+B_W:=h^*(K_X+B)$. Since $\Nklt(X,B)=\Nklt(X,B,\Mm)$, there exists $0<\delta\ll 1$ such that $(W,B_W+\delta E)$ is sub-lc and $\frac{1}{2}h^*A-\delta E$ is ample$/U$. Thus $\Mm_W+\frac{1}{2}h^*A-\delta E$ is ample$/U$, and we may pick $0\leq H_W\sim_{\Rr,U}\Mm_W+\frac{1}{2}h^*A-\delta E$ such that $(W,B_W+H_W+\delta E)$ is sub-lc. Let $\Delta':=h_*(B_W+H_W+\delta E)$. Then $(X,\Delta')$ is lc and $\Delta'\sim_{\Rr,U}B+\Mm_X+\frac{1}{2}A$. Possibly replacing $A$ we may assume that $(X,\Delta'+\frac{1}{2}A)$ is lc. By \cite[Theorem 1.5]{HH20},  $(X,\Delta'+\frac{1}{2}A)$ has a good minimal model. By \cite[Lemma 4.2]{HL21a}, we get (2).
\end{proof}

\begin{proof}[Proof of Theorem \ref{thm: bir12 1.1 gpair}]
If $K_X+B+\Mm_X$ is not pseudo-effective$/U$, then the theorem follows from \cite[Lemma 4.4(1)]{BZ16} after passing to a dlt model of $(X,B,\Mm)$. So we may assume that $K_X+B+\Mm_X$ is pseudo-effective$/Z$. By \cite[Theorem 3.14]{HL21a}, we only need to prove (2), so we may assume that $(X,B,\Mm)$ is $\Qq$-factorial dlt. We run a $(K_X+B+\Mm_X)$-MMP$/U$ with scaling of an ample$/U$ $\Rr$-divisor $H\geq 0$:
$$(X,B,\Mm):=(X_0,B_0,\Mm)\dashrightarrow (X_1,B_1,\Mm)\dashrightarrow\dots\dashrightarrow (X_i,B_i,\Mm)\dashrightarrow\dots.$$
By Theorem \ref{thm: Has22b 4.1 rel ver} ($U$ and $Z$ in Theorem \ref{thm: Has22b 4.1 rel ver} both correspond to our $U$, $A_Z$ and $A$ of Theorem \ref{thm: Has22b 4.1 rel ver} correspond to $0$, and $C$ corresponds to our $A$), 
$K_{X_i}+B_i+\Mm_{X_i}$ is log abundant$/U$ with respect to $(X_i,B_i,\Mm)$ for every $i$. By \cite[Theorem 2.24]{HL21a} and Theorem \ref{thm: Has22a 3.15 rel ver}, this MMP terminates with a log minimal model of $(X,B,\Mm)/U$.
\end{proof}

\end{document}